\pdfoutput=1
\RequirePackage{ifpdf}
\ifpdf % We~are running pdfTeX in pdf mode
\documentclass[pdftex]{sigma}
\else
\documentclass{sigma}
\fi

\usepackage{mathtools}
\usepackage{tikz}
\usepackage{tikz-cd}
\usetikzlibrary{positioning, automata, patterns, decorations.pathreplacing, arrows, shapes, calc}

\usepackage[capitalize]{cleveref}

\usepackage{macros} % package with our theorem styles and math commands

\newtheorem{Theorem}{Theorem}[section]
\newtheorem*{Theorem*}{Theorem}
\newtheorem{Corollary}[Theorem]{Corollary}
\newtheorem{Lemma}[Theorem]{Lemma}
\newtheorem{Proposition}[Theorem]{Proposition}
\newtheorem{Conjecture}[Theorem]{Conjecture}
 { \theoremstyle{definition}
\newtheorem{Definition}[Theorem]{Definition}

\newtheorem{Example}[Theorem]{Example}
\newtheorem{Remark}[Theorem]{Remark}
\newtheorem{Question}[Theorem]{Question}
}

\begin{document}

\allowdisplaybreaks

\newcommand{\arXivNumber}{2304.00699}

\renewcommand{\PaperNumber}{073}

\FirstPageHeading

\ShortArticleName{$\widehat{Z}$ and Splice Diagrams}

\ArticleName{$\boldsymbol{\widehat{Z}}$ and Splice Diagrams}

\Author{Sergei GUKOV~$^{\rm a}$, Ludmil KATZARKOV~$^{\rm b}$ and Josef SVOBODA~$^{\rm a}$}

\AuthorNameForHeading{S.~Gukov, L.~Katzarkov and J.~Svoboda}

\Address{$^{\rm a)}$~Department of Mathematics, California Institute of Technology, Pasadena, CA 91125, USA}
\EmailD{\href{mailto:gukov@caltech.edu}{gukov@caltech.edu}, \href{mailto:svo@caltech.edu}{svo@caltech.edu}}

\Address{$^{\rm b)}$~Department of Mathematics, University of Miami, Coral Gables, FL 33146, USA}
\EmailD{\href{mailto:l.katzarkov@miami.edu}{l.katzarkov@miami.edu}}

\ArticleDates{Received November 22, 2024, in final form August 16, 2025; Published online August 26, 2025}

\Abstract{We study quantum $q$-series invariants of 3-manifolds $\widehat{Z}_\sigma$ of Gukov--Pei--Putrov--Vafa, using techniques from the theory of normal surface singularities such as splice diagrams. We show that the (suitably normalized) sum of all $\widehat{Z}_\sigma$ depends only on the splice diagram, and in particular, it agrees for manifolds with the same universal abelian cover. We use these ideas to find simple formulas for $\widehat{Z}_\sigma$ invariants of Seifert manifolds. Applications include a better understanding of the vanishing of the $q$-series $\widehat{Z}_\sigma$. Additionally, we study moduli spaces of flat $\operatorname{SL}_2(\mathbb{C})$ connections on Seifert manifolds and their relation to spectra of surface singularities, extending a result of Boden and Curtis for Brieskorn spheres to Seifert rational homology spheres with 3 singular fibers and to Seifert homology spheres with any number of fibers.}

\Keywords{3-manifold topology; quantum invariant; surface singularity; splice diagram}

\Classification{57K31; 32S50; 32S25; 32S55}

\section{Introduction}\label{sec:intro}

Low-dimensional topology and singularity theory have been always closely related, with singularities of curves and surfaces providing interesting examples of knots and 3-manifolds in the first place. More recently, the techniques from singularity theory have been very fruitful in the study of invariants coming from gauge theory, such as Seiberg--Witten invariants and Heegaard--Floer homology of Ozsvath and Szab\'o, in the work of N\'emethi and Nicolaescu \cite{NeNi02}. In this paper, we~apply those techniques to \emph{quantum} invariants of 3-manifolds.

More specifically, we study $q$-series invariants $\Zhat_{\sigma}(Y; q)$ of a 3-manifold $Y$, the so-called GPPV invariants, associated with quantum groups $\mathcal{U}_q (\mathfrak{g})$ at generic $|q|<1$ and labeled by a \spc{} structure $\sigma \in \spc(Y)$ on $Y$. They are expected to have a categorification \cite{GPV16}.

Mostly interested in topological and geometric aspects, throughout the paper, we consider the simplest non-trivial choice of $\mathfrak{g} = \mathfrak{sl}_2$, which corresponds to the ``gauge group'' $G = \operatorname{SU}(2)$, though much of the discussion can be generalized to higher rank root systems.\footnote{We note that by ``gauge group'' one can mean at least two different things. One is the gauge group in complex Chern--Simons theory, for which $\Zhat_{\sigma}(Y; q)$ provides a non-perturbative definition that behaves well under cutting and gluing. In contrast, the other one is a gauge group in 3d-3d correspondence and string theory fivebrane setup. The former group is a complexification of the latter, and it is the latter that we use to label $\Zhat_{\sigma}(Y; q)$ in this paper. So, even when we talk about $G = \operatorname{SU}(2)$, the gauge group in complex Chern--Simons theory and its relations to quantum groups is, in fact, $\operatorname{SL}_2(\C)$.} Also, unless stated otherwise, throughout the paper, we work in the category of connected irreducible oriented negative definite plumbed 3-manifolds \cite{Neu81} with $b_1 = 0$, for which one can use the definition of~$\Zhat_{\sigma}(Y; q)$ as in~\cite{GPPV20}. (A reader interested in various generalizations and extensions is welcome to consult \cite{AJK21,Chae21, GM21, Pa20a, Pa20b, Pa21, Ri22}.)

Our main goal is to approach the invariants $\Zhat_{\sigma}(Y; q)$ by using various methods of the singularity theory of normal complex surfaces, most notably the universal abelian covers and splice diagrams. Among other things, this new perspective leads to simpler expressions for $\Zhat_{\sigma}(Y; q)$ for certain families of 3-manifolds (see Theorem~\ref{thm:seifert}) and sheds new light on what topological information they capture.
Specifically, $\Zhat_{\sigma}(Y; q)$ have close cousins $Z_a(Y; q)$ (without a hat!) that are not expected to admit a categorification, but are nevertheless more natural from the viewpoint of complex Chern--Simons theory \cite{GMP16}. While both sets of invariants can be defined independently, they are linearly related and, therefore, can be expected to contain roughly the same topological information about $Y$. By taking a closer look at the structure of the linear relation between $\Zhat_{\sigma}(Y; q)$ and $Z_a(Y; q)$ through the looking glass of the singularity theory, we~observe that $Z_0(Y; q)$ is, surprisingly, a much simpler invariant of $Y$ in the precise sense that we explain below.

We hope that a fresh new look through the lens of the singularity theory in the future will lead to a better understanding of the structure of the invariants $\Zhat_{\sigma}(Y; q)$, just as it led to simplified formulae for Seifert manifolds in this work. The fact that singularity theory involves complex surfaces brings us one step closer to one of the main motivations in studying these invariants, namely developing new homological 3-manifold invariants that could help us explore the world of smooth 4-manifolds.

\subsection{Main results}\label{ss:main_results}

Let $Y$ be a closed, oriented, irreducible 3-manifold which is plumbed, see Section~\ref{sec:plumbing} for definitions. Moreover, assume that the plumbing matrix can be chosen to be negative definite. We call such a 3-manifold \emph{a negative definite plumbed manifold} as is traditional.\footnote{The reader should keep in mind that negative definiteness is a property of the plumbing graph rather than of~$Y$.} Note that these conditions are equivalent to $Y$ being a link of an isolated normal surface singularity $X$. Moreover, let us assume that $Y$ is a rational homology sphere, i.e., $b_1(Y)=0$, so it has finite $H:=H_1(Y,\ZZ)$. Under these assumptions, $\Zhat_\sigma(q)$ invariants were defined in \cite{GPPV20}.

The \emph{splice diagram} is a combinatorial notion developed by Neumann, Eisenbud and Wahl \mbox{\cite{EN86, NW05a}} for the study of singularities and 3-manifolds described above. Constructed from the plumbing graph, it contains certain essential part of it and proves very useful for studying of $\Zhat_\sigma(q)$ invariants. The following theorem states that the series $Z_0(q) = \sum_{\sigma \in \spc(Y)}{\Zhat_\sigma(q)}$ essentially depends on the splice diagram (see Section~\ref{ssec:splice_diagrams} for definitions). This applies in particular to homology spheres, where it provides a useful computational tool for obtaining the unique~${\Zhat_\sigma(q) = Z_0(q)}$. Denote by $\lambda(Y)$ the Casson--Walker invariant \cite{Wa92}.

\begin{Theorem}
	Let $Y$ be a $3$-manifold satisfying the assumptions above. The $q$-series
	\begin{equation}
		q^{-6 \lambda(Y)} Z_{0}\bigl(q^{|H|}\bigr)
	\end{equation}
	depends only on the splice diagram of $Y$.
\end{Theorem}

Topologically, splice diagrams determine universal abelian covers and vice versa \cite{NW05a, Pe09}. \emph{Universal abelian cover} of a manifold $Y$ is the maximal cover of $Y$ with abelian deck group and by our assumptions, this cover is finite. We obtain the following corollary of the theorem.

\begin{Corollary}
	If $Y_1$ and $Y_2$ have the same universal abelian cover, then
	\[
		q^{-6 \lambda(Y_1)} Z_0\bigl(Y_1, q^{|H_1(Y_1)|}\bigr) =
		q^{-6 \lambda(Y_2)} Z_0\bigl(Y_2, q^{|H_1(Y_2)|}\bigr).
	\]
\end{Corollary}
For example, Seifert manifolds $Y_1 = M(1;(7,1),(7,1),(7,4))$, $Y_2 = M(1;(7,1),(7,2),(7,3))$ satisfy $Y_1^{\rm ab}=Y_2^{\rm ab}$ and $H=(\ZZ/7\ZZ)^2$. Then $\lambda(Y_1)=0$, $\lambda(Y_2)=-21/2$ (see \eqref{eq:casson_formula}) and
\[
	Z_0\bigl(Y_1,q^{49}\bigr) = q^{217/2} \bigl(1 - 3 q^{35} + 3 q^{84} - q^{147} + q^{539} +\cdots\bigr) = q^{63} Z_0\bigl(Y_2,q^{49}\bigr).
\]

Seifert manifolds (Seifert fiber spaces) mentioned above are circle fibrations over 2-dimen\-sion\-al orbifolds. They admit ``star-shaped'' plumbing graphs and form basic building blocks of plumbed manifolds in Jaco--Shalen--Johannson decomposition. In this special case, we can go further with simplifications of $Z_0(Y)$ and we will give an explicit and rather simple formula, the ``reduction'' Theorem~\ref{thm:seifert}, using the data of Seifert fibration (rather than that of the plumbing graph). Using the group action of $H$ on $Y^{\rm ab}$, we will extract the separate $\Zhat_{\sigma}(q)$ for each
$\sigma$ in~$\spc(Y)$.

The next theorem gives explicit formulas for $\Zhat_\sigma(q)$ invariants of Seifert manifolds, see Section~\ref{sec:seifert} for definitions and conventions. We introduce the ``Laplace transform'' $\LC_A$, defined on~monomials as
\smash{$
	\LC_A(t^n) = q^{\frac{n^2}{4A}}
$}
and extended linearly to any formal power series in $t$. For a~rational function $f(t)$, we denote $\se f(t)$ the ``symmetric expansion'', which is the average of the expansions of the given rational function as $t \ra 0$ (as a Laurent power series in $t$) and as $t \ra \infty$ \big(as a Laurent power series in $t^{-1}$\big). Finally, we pick a reference ``anticanonical'' $\spc{}$ structure~$\ac$ to express all $\Zhat_\sigma(q)$ using the natural action of $H = H_1(Y,\ZZ) = \langle g_0,\dots, g_k \rangle$ on $\spc(Y)$ (see Sections \ref{spc_structures} and \ref{sec:seifert} for notations and definitions).

\begin{Theorem}[reduction theorem]
	Let $Y = M(b;(a_1, b_1), (a_2, b_2), \dots, (a_k, b_k))$ be a Seifert manifold over $S^2$ $(k \geq 3)$ with
\[
e=-b+\sum_{i=1}^k b_i/a_i<0.
\]
 We set $A=a_1 a_2\cdots a_k$ and $A_i=A/a_i$ for $i=1,\dots, k$. Then
	\begin{gather}\label{eq:red_theorem_z0}
		Z_0\bigl(q^{|H|}\bigr) = q^\Lambda \LC_A(\se f_0(t)),
	\end{gather}
	where
	\[
		f_0(t) = \frac{\bigl(t^{A_1}-t^{-A_1}\bigr)\bigl(t^{A_2}-t^{-A_2}\bigr)\cdots\bigl(t^{A_k}-t^{-A_k}\bigr)}{\bigl(t^A-t^{-A}\bigr)^{k-2}},
	\]
 and $\Lambda = \Lambda(Y)$ is an explicit rational number given in \eqref{eq:Lambda}.
	Moreover, for $h \in H$ we have
	\[
		\Zhat_{h \ac}(q) = q^{\Lambda/|H|} \LC_{A|H|}(\se f_{h \ac}(t)),
	\]
	where
	\[
		\sum_{h \in H} f_{h \ac}(t)h = t^{(k-2)A - \sum_i A_i}\frac{\bigl(g_1 t^{2A_1}-1\bigr)\cdots\bigl(g_k t^{2A_k}-1\bigr)}{\bigl(g_0 t^{2A}-1\bigr)^{k-2}}.
	\]
\end{Theorem}

These formulae have both computational and conceptual significance. While in the original formulation of $\Zhat_{\sigma}(q)$ we need to work with the plumbing data, which can be very large due to the presence of continued fractions, here we get $\Zhat_{\sigma}(q)$ simply from the data of the Seifert fibration. It also sheds light on the role of the action of $H$, leading to a better understanding of the vanishing of the $q$-series for certain $\sigma$ in Corollary~\ref{cor:vanishing}. Note that the independence of $f_0(t)$ on $b_i$ in \eqref{eq:red_theorem_z0} is an instance of Theorem~\ref{thm:plumbing}.

\subsection{Spectra in algebra and geometry}

Finally, in Section~\ref{sec:spectrum}, we tie together several aspects of the story, mainly focusing on Seifert manifolds. These involve geometric structures, such as flat connections on $Y$ and the invariants of the corresponding Brieskorn-type complete intersection singularity, as well as vertex algebras.

The $\Zhat_{\sigma}(Y)$ invariants are defined and studied using plumbing graphs. In singularity theory, negative definite plumbing graphs correspond to resolutions of the singularity which are smooth complex surfaces. Many deep results in singularity theory concern other smooth surfaces related to the singularity, namely those which are results of a smoothing. The theory of Milnor fibrations was originally developed for singularities of a single holomorphic function (hypersurface singularities) \cite{Mil69} and later extended to complete intersections \cite{Loo84}.

\emph{Spectrum of hypersurface singularity} \cite{St85} was defined by Arnold, Steenbrink and Varchenko and generalized later by Steenbrink and Ebeling to complete intersections \cite{EbSt98}. The spectrum is a collection of real numbers that refine the eigenvalues of the monodromy, using a natural mixed Hodge structure on the cohomology of Milnor fiber. Its most important property is the upper-semicontinuity under deformations.

Recall that for pairwise coprime integers $p$, $q$, $r$, the Brieskorn homology sphere $Y = \Sigma(p, q, r)$ is the link of the Brieskorn singularity $X\colon x^p+y^q+z^r = 0$. A relation between the Milnor fiber of $X$ and the topology of the link was given by Fintushel and Stern \cite{FS90}.

\begin{Theorem}[{\cite[Theorem 2.10]{FS90}}]\label{thm:FS}
	Let $\lambda(Y)$ be the $\operatorname{SU}(2)$ Casson invariant of $Y$ and $\sigma(X)$ the signature of the Milnor fiber of $X$. Then
$\lambda(Y) = \frac{\sigma(X)}{8}$.
\end{Theorem}
This relation, the celebrated `Casson invariant conjecture' \cite{NeuWa90} has been extended in several directions, e.g., to the singularities of splice-type \cite{NeOk08}. In a similar spirit, Curtis \cite{Cu01} defined the~$\operatorname{SL}_2(\C)$ Casson invariant \smash{$\lambda^{C}_{\operatorname{SL}_2(\C)}(Y)$} and with Boden proved a relation between \smash{$\lambda^{C}_{\operatorname{SL}_2(\C)}(Y)$} and the Milnor number $\mu(X)$, the rank of the cohomology of the Milnor fiber~\cite{BC06}.

\begin{Theorem}[{\cite[Theorem~2.4.]{BC06}}]\label{CassonMilnorIntro}
With the notation above,
\smash{$
\lambda^C_{\operatorname{SL}_2(\C)}(Y) = \frac{\mu(X)}{4}$}.
\end{Theorem}

We extend this relation in two directions -- to Seifert rational homology spheres with 3 singular fibers, and to Seifert homology spheres with any number of singular fibers.

For a Seifert rational homology sphere $Y = M(b;(a_1, b_1), \dots, (a_3,b_3))$, we need to consider the spectrum of the universal abelian cover singularity $X^{\rm ab}$, given by the equation
$
 X^{\rm ab}\colon x_1^{a_1}+x_2^{a_2}+x_3^{a_3}=0$,
where $a_1$, $a_2$, $a_3$ are not necessarily coprime. Instead of the Milnor number, we consider the number of orbits of a natural non-free $(\Z_2)^2$-action of the spectrum $\spec\bigl(X^{\rm ab}\bigr)$ of~$X^{\rm ab}$ (see \eqref{spec_three} and \eqref{eq:action}). The action replaces the factor of $1/4$ in Theorem~\ref{CassonMilnorIntro}.

On the topology side, we need the number \smash{$\lambda^{\rm nab}_{\operatorname{SL}_2(\C)}(Y)$} of characters of non-abelian flat connections, including the reducible ones (not present for Seifert homology spheres). This number was computed by Cui, Qiu and Wang \cite{CuQiWa21}.

\begin{Proposition}[Proposition~\ref{prop:rat_sph_spec}]
Let $Y = M(b;(a_1, b_1), \dots, (a_3,b_3))$ be a Seifert rational homology sphere. With the notation above, we have
\[
	\lambda^{\rm nab}_{\operatorname{SL}_2(\C)}(Y) =\left\lvert\frac{\spec\bigl(X^{\rm ab}\bigr)}{(\Z_2)^2}\right\rvert.
\]
\end{Proposition}

For Seifert homology spheres, we use an alternative definition of Casson invariant \smash{$\lambda^{P}_{\operatorname{SL}_2(\C)}$} by Abouzaid and Manolescu \cite{AbMa20}, which takes into account higher-dimensional components of the moduli space of $\operatorname{SL}_2(\C)$ flat connections on $Y$. Brieskorn singularity is replaced by the Brieskorn--Hamm complete intersection (see \eqref{icis}).
	\begin{Theorem}[Theorem~\ref{thm:seifert_casson}]
		Let $Y$ be a Seifert homology sphere and $X$ the corresponding Brieskorn-Hamm complete intersection singularity with Milnor number $\mu(X)$. Then
		\begin{equation}\label{casson_milnor_intro}
			\lambda^{P}_{\operatorname{SL}_2(\C)}(Y) = \frac{\mu(X)}{4}.	
		\end{equation}
	\end{Theorem}
The proof is merely a combination of formulas by Hamm \cite{Ham72} for Milnor number and Boden and Yokogawa \cite{BoYo96} for the flat connections, but we could not find it in the literature.\footnote{Boden and Curtis \cite[p.~10]{BC06} related the Euler characteristic of the union of the top-dimensional components of the moduli space with the Milnor number of a related \emph{hypersurface} singularity.}$^{,}$\footnote{Recently, Mu\~noz derived the same result using gauge theory \cite[Corollary B]{munoz25}.}

If one asks what is the topological analogue of the spectrum or the eigenvalues of the monodromy, the closest seem to be the ``rotation numbers'' of the components of the moduli space of connections. In the course of the proof of Theorem~\ref{thm:seifert_casson}, it is quite amusing to notice that not only these numbers match, but also the multiplicities of the spectral numbers match the Euler characteristic of the higher-dimensional components of the moduli space of connections. This suggests a deeper structure behind the numerical relation \eqref{casson_milnor_intro}.

The above relation between the monodromy of a singular surface with the moduli spaces of representations suggests the existence of new nonabelian noncommutative Hodge structure --- a~combination of Higgs bundles and Landau--Ginzburg theory, and related to it a new, derived spectrum. We sketch some examples and pose some questions about these structures. This more Hodge theoretic, rather than topological construction opens new ways of studying degenerations of nonabelian Hodge structures and a potential Landau--Ginzburg theory interpretation of~$\Zhat_{\sigma}(q)$.

On the algebra side, namely in vertex algebra (VOA), there is a different notion of the spectrum. It also refers to a collection of rational numbers, $\{ \Delta_i \}$, that determine leading $q$-powers in the $q$-expansion of VOA characters. These numbers are called {\it conformal weights} or {\it conformal dimensions} since they are defined as eigenvalues of the conformal vector $L_0$, which is part of the mathematical definition of vertex algebra.

Since one of the predictions of the so-called 3d-3d correspondence is that a closed 3-manifold $Y$ corresponds to a VOA \cite{CCFGH18, CCFFGHP22, Sug22} for which $\Zhat_{\sigma}$ is a character, it is natural to ask how the spectrum of conformal weights $\Delta_{\sigma}$ relates to other spectra mentioned above, defined more geometrically. Contrary to what one might naively expect, we find that the spectrum of conformal weights~$\Delta_{\sigma}$ is rather different from the spectrum of the corresponding hypersurface singularity: while the former determines the leading $q$-power in the $q$-expansion, the latter encodes the structure of the $q$-series coefficients at large $q$-powers. In turn, the latter determines the behavior of the $q$-series near $q = {\rm e}^{\hbar} \approx 1$ or, equivalently, the expansion in $\hbar$ near $\hbar = 0$.

To summarize, expansion in $q$, as in $\Zhat_{\sigma} (q)$, has direct contact with counting problems (curve counting, BPS states, etc.), with vertex algebras, and with the spectrum of conformal dimensions. These connections become less natural near $q = {\rm e}^{\hbar} \approx 1$, where connections to complex Chern--Simons theory and the spectrum of hypersurface singularities become manifest. The two expansion limits are related by resurgent analysis, a powerful technique that, roughly speaking, allows to transfer (enlarge) the domain of a given function (or a power series, possibly, with zero radius of convergence).

This interplay between expansions near $q = 0$ and $q = {\rm e}^{\hbar} \approx 1$ plays an important role in curve counting \cite{DT, GV, Katz, MNOP} and --- since $\Zhat_{\sigma} (q)$ admit an interpretation via curve counting too \cite{EGGKPSS22, EGGKPS22, GPPV20, GPV16} --- in the study of $\Zhat_{\sigma}(q)$-invariants \cite{AM22, Chun17, Chung20, GMP16,Wu20}. In particular, near $q = {\rm e}^{\hbar} \approx 1$ more natural objects are $Z_a (q)$, without a `hat', and trans-series $Z_{\alpha} (q)$ of the complex Chern--Simons theory. This is consistent with what we find here: $Z_0$ and, more generally, $Z_{\alpha}$ appear to be more natural from the perspective of the singularity theory. And, perhaps not surprisingly (given the above explanations), the spectrum of the hypersurface singularity associated with $Y$ is closely related to the set of values of the classical Chern--Simons functional on $Y$, i.e., to the singularities on the Borel plane.

\section{Plumbed manifolds}\label{sec:plumbing}

In this section, we study $\Zhat_{\sigma} (q)$ invariants and their modifications. We use the definition given by the integral formula \cite{GM21,GPPV20} and interpret its parts, building analogies with invariants that appear in the study of normal surface singularities, such as the topological Poincar\'e series.\looseness=1

\begin{figure}
	\centering
	\includegraphics[scale = 1.3]{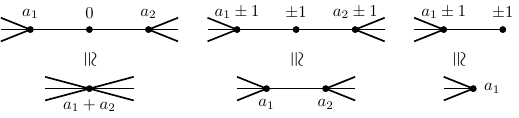}
	\caption {Moves on plumbing trees that preserve the $3$-manifold.}
	\label{fig:kirby}
\end{figure}

Recall the theory of plumbed 3-manifolds. Let $\Gamma$ be a finite tree (a graph with no cycles).
We call \emph{leaf} a vertex of degree one and \emph{node} a vertex of degree 3 and more. Additionally, we assign an integer label (framing) $m_v$ to each vertex $v$. Associated with this data, there is a plumbed manifold $Y$: For each vertex $v$, we take a circle bundle over 2-sphere with Euler number $m_v$ and then we glue along tori corresponding with the edges. One can think of $Y$ as the boundary of a 4-manifold $X$ constructed by gluing disc bundles in a similar fashion. It is convenient to describe $\Gamma$ by the \emph{plumbing matrix} $M$
\[
	M_{vw} =
	\begin{cases}
		m_v & \text{if} \ v = w , \\
		1 & \text{if} \ (v, w) \in \text{Edges}(\Gamma), \\
		0 & \text{otherwise}.
	\end{cases}
\]
Two different plumbing graphs represent the same 3-manifold $Y$ if and only if they are related to each other by Kirby--Neumann moves \cite{Neu81}; see Figure~\ref{fig:kirby}.

The main source of plumbed manifolds are the isolated normal singularities of complex surfaces. To any such singularity germ $(X, o)$, we can assign a 3-manifold (``the \emph{link} of the singularity'') by embedding a neighborhood of the singular point into $(\CC^n, 0)$ and intersecting with a sufficiently small real sphere around zero (see, e.g., \cite{Nem22} for details of the construction). Such plumbing graphs are always negative definite in the sense that the plumbing matrix is. Conversely, every negative definite plumbing graph is a resolution graph of some singularity and the corresponding 3-manifold is its link.

\begin{Example}
	The singularity $E_{12}\colon x^2+y^3+z^7 = 0$ has a resolution with the intersection matrix
	\[
		\begin{pmatrix}
				-7 & 0 & 0 & 1 \\
				0 & -3 & 0 & 1 \\
				0 & 0 & -2 & 1 \\
				1 & 1 & 1 & -1
			\end{pmatrix}.
	\]
	This matrix is correspondingly a negative definite plumbing matrix of the link of $E_{12}$, Brieskorn sphere $\Sigma(2, 3, 7)$.
\end{Example}

\subsection[Lattices and spinc structures]{Lattices and $\boldsymbol{\spc{}}$ structures}\label{spc_structures}

Before we get into the $\Zhat_{\sigma} (q)$ invariants, we need some preparation on lattices associated with plumbed manifold $Y$ and \spc{} structures on $Y$. We will use $\sigma$, $\sigma'$ to denote \spc{} structures and $g, h, \dots$ for elements of $H = H_1(Y, \ZZ)$. Denote $\hat{H} = \operatorname{Hom}(H, \C^*)$ the Pontryagin dual of~$H$. For a rational homology sphere $Y$, $H$ is finite and acts freely and transitively on the set of \spc{} structures $\spc(Y)$, so $|\spc(Y)| = |H|$. We denote the action of $h$ on $\sigma$ by~$h\sigma$.

Assume now that $\tilde{X}$ is a plumbed 4-manifold with boundary $Y$ and plumbing matrix $M$ of size $s \times s$. We consider the lattice $L = H_2\bigl(\tilde{X}, \ZZ\bigr)$ equipped with the intersection form. This form induces an embedding of $L$ to $L' = H^2\bigl(\tilde{X}, \ZZ\bigr)$ and extends naturally to $L'$. We identify the lattice $L'$ with \smash{$\bigl(\ZZ^s, \bigl(\vl, M^{-1} \vl'\bigr)\bigr)$} and $L$ with $M\ZZ^s \subset L'$.\footnote{The basis $E^*_i$, $i = 1, \dots, s$ (resp.\ $E_i$) of $L'$ (resp.\ $L$) used in~\cite{Nem22} corresponds to the vectors $-e_i$ (resp.\ $M e_i$) where $e_i$ is the standard basis of~$\ZZ^n$.}
The group $H$ is then naturally identified with $L'/L = \ZZ^s/M \ZZ^s$ using the short exact sequence
\[
	0 \rightarrow H^2\bigl(\tilde{X}, Y, \ZZ\bigr) \rightarrow H^2\bigl(\tilde{X}, \ZZ\bigr) \rightarrow H^2(Y, \ZZ) \rightarrow 0
\]
and Poincar\'e--Alexander duality in the first and third term. As $M$ is negative definite, the order of $H$ is $|H|=\det(-M)$.

Let $\vm$ and $\vdelta$ be the vectors of the framings and the degrees of the vertices, respectively. The \spc{} structures on $Y$ are naturally identified with the \emph{characteristic vectors}, which are the elements of the set
$\Char(Y) = (2\ZZ^s+\vm)/2M\ZZ^s$.
As in \cite{GM21}, we use the identity $\vdelta + \vm = M\vu,$ where $\vu = (1, 1, \dots, 1)$,
to further identify $\Char(Y)$ with a slightly different set $\Char'(Y)$
\begin{equation}\label{eq:char}
	\Char'(Y) = \bigl(2\ZZ^s+\vdelta\bigr)/2M\ZZ^s,
\end{equation}
via the map induced by $\vl \mapsto \vl - M\vu$.

The set of \spc{} structures has natural involution $\sigma \rightarrow \bar \sigma$ which on characteristic vectors (either $\Char$ or $\Char'$) acts simply by $\vsigma \mapsto -\vsigma$. Fixed points of this involution, i.e., the intersection~${(2\ZZ^s+\delta) \cap M\ZZ^s}$, are naturally identified with \spin{} structures. %chern class

In the context of links of singularities, a special choice of \spc{} structure plays an important role \cite{NeNi02}. This ``canonical'' \spc{} structure \can{} has characteristic vectors $\vm + 2\vu \in \Char$ and~${2\vu - \vdelta}$ in $\Char'$. A more suitable choice for our purposes is the ``anticanonical'' \spc{} structure $\ac := \overline{\can{}},$ with characteristic vector $\vdelta - 2\vu$ in $\Char'$.

The choice of \spc{} structure $\ac$ gives us a $H$-equivariant map $H \rightarrow \spc(Y)$ given by $h \rightarrow h\ac$. In terms of a vector representative $\vh$ of $h$ in $H \cong \ZZ^s/M\ZZ^s$ and $\vdelta - 2\vu \in \Char'(Y)$, the action reads $\vh \mapsto \vdelta - 2\vu + 2\vh$. Note that this map preserves the involutions if and only if~$\ac$ is a \spin{} structure. This is not the case in general.

\subsection[The q-series invariants]{The $\boldsymbol{q}$-series invariants}\label{Zhat_review}

Following \cite{GM21}, we fix a plumbing tree $\Gamma$ with $s$ vertices, denote vectors in $\vl \in \ZZ^s$ by letters with arrows and write $l_v$ for respective component of a vertex $v$.

Let $Y$ be a rational homology sphere given by a negative definite plumbing\footnote{In \cite{GM21}, a definition was given for a more general class of plumbing graphs called weakly negative plumbings. However, it can be shown that every weakly negative graph can be transformed to a negative definite graph by a~sequence of Kirby--Neumann moves, by a~modification of the argument used for the characterization of negative definite plumbings by Eisenbud and Neumann \cite{EN86}. Therefore, we use definite graphs only.} with a fixed plumbing graph $\Gamma$ on $s$ vertices and plumbing matrix $M \in \GL(\ZZ, s)$. Denote by $\vdelta$ the vector in~$\ZZ^s$ of degrees of vertices of $\Gamma$.

The $q$-series $\Zhat_{\sigma} (q)$ are defined as\footnote{Let us note that by the Weyl $\ZZ_2$ symmetry of the integrand, we have $\Zhat_{\sigma} = \Zhat_{\bar{\sigma}}$ and they are often identified. We do not do it here, i.e., use the so-called ``unfolded'' version \cite{CCFGH18}.}
\begin{equation}\label{eq:plumbing1}
	\Zhat_{\sigma} (q) = q^{\frac{-3s-\Tr(M)}{4}} \cdot \text{v.p.} \oint\limits_{|z_v| = 1} \prod_{v\in \Ver} \frac{{\rm d}z_v}{2\pi{\rm i} z_v}\left( z_v - \frac{1}{z_v} \right)^{2-\deg(v)} \cdot \Theta_{\sigma}^{-M}(q, \vz),
\end{equation}
where
\[%\label{Thetaa}
	\Theta_{\sigma}^{-M}(q, \vz) = \sum_{\vl \in 2M\Z^s + \vsigma} q^{-\frac{(\vl, M^{-1}\vl)}{4}} \prod_{ v\in \Ver} z_v^{l_v}.
\]

Here, v.p.\ denotes taking the principal value of the integral, given by the average of the
integrals over the circles $|z_v| = 1 + \epsilon$ and $|z_v| = 1 - \epsilon$, for $\epsilon > 0$ small. The negativity condition is needed for the convergence of this $q$-series. The label $\sigma$ stands for a choice of {\spc} structure on $Y$, identified with a vector in $\Char'(Y) = \bigl(2\ZZ^s+\vdelta\,\bigr)/2M\ZZ^s$; see \eqref{eq:char}.

We can reformulate the integral formula using symmetric expansions. Recall that for a~rational function
$r(t) \in \C(t)$, the \emph{symmetric expansion} of $r(t)$, denoted $\se r(t)$, is the average of Laurent expansions around $0$ and $\infty$
\[
	\se r(t) := \frac12 \bigl(\expn_{t\rightarrow 0} r(t) + \expn_{t \rightarrow \infty} r(t)\bigr).
\]
Using this notion, we can turn \eqref{eq:plumbing1} into the following formula:
\[%\label{eq:plumbing2}
	\Zhat_\sigma(q) =
	2^{-s} q^{\frac{-3s-\Tr(M)}{4}} \LC_{-M}(F_\sigma(\vz )),
\]
where $F_\sigma(\vz ) = F_\sigma(z_1,\dots,z_s) \in \ZZ\big[\big[z_1^{\pm 1}, \dots, z_s^{\pm 1}\big]\big]$ are formal power series
\begin{equation}\label{eq:zhat_expansion}
	F_{\sigma}(\vz )=
	\pi_\sigma \left( \prod_{v \in \Ver} 2 \se \left(z_v-\frac{1}{z_v}\right)^{2-\delta_v} \right),
\end{equation}
where $\pi_\sigma$ is the projection on the class of $\sigma$
\[
	\pi_\sigma \bigl(\vz\,{}^{\vl}\bigr) =
	\begin{cases}
		\vz\,{} ^{\vl} & \text{if } [\vl ] = \sigma, \\
		0 & \text{otherwise}.
	\end{cases}
\]
$\LC_{-M}$ denotes the following transformation (``Laplace transform'')
\smash{$\LC_{-M}\colon \vz\,{}^{ \vl} \rightarrow q^{\frac{-(\vl, M^{-1}\vl)}{4}}$},
extended linearly to formal power series in $\vz\,{}^{\pm 1}$.

We can eliminate the projections using the \spc{} structure \ac
\begin{equation}\label{eq:zhat_expansion_can}
	\sum_{h \in H} F_{h\ac}(\vz ) h=
	 \prod_{v \in \Ver} z_v^{\delta_v - 2} \prod_{v \in \Ver}
	2\se \bigl(g_v z_v^2-1\bigr)^{2-\delta_v},
\end{equation}
where $g_v$ are the generators of $H$ given by classes $[e_v]$ of canonical basis vectors of $\ZZ^s/M\ZZ^s$.
Both sides of this equation should be understood as elements of the group ring $\ZZ[\vz\,{}^{ \pm 1} ][H]$ of $H,$ with coefficients in formal power series in $\vz$ and the symmetric expansion is naturally extended to $\CC(t)[H]$.

\section[Z\_0 and splice diagrams]{$\boldsymbol{Z_0}$ and splice diagrams}\label{sec:splice}

In this section, we will concentrate on the series $Z_0(q)$ defined as the sum of $\Zhat_\sigma(q)$ over all $\sigma \in \spc(Y)$,
$Z_0(q) = \sum_{\sigma \in \spc(Y)} \Zhat_\sigma(q)$.
Clearly, one has
\[
	Z_0(q) =
	2^{-s} q^{\frac{-3s-\Tr(M)}{4}} \LC_{-M}(F(\vz )),
\]
where $F(\vz )$ is given by the expansion \eqref{eq:zhat_expansion} without the projection $\pi_\sigma$.

We will prove that $Z_0(q)$ can be reconstructed\footnote{Up to an overall power.} from less data than the full plumbing graph (see Theorem~\ref{thm:plumbing}). The relevant notion is the splice diagram \cite{EN86,Neu81,Sie80}, which was originally developed for the study of singularities but proves to be very useful in our context as well. For homology spheres, the plumbing graph can be reconstructed from the splice diagram, but the splice diagram can be used for faster computations of the unique $\Zhat_\sigma(q) = Z_0(q)$ in this case.

\subsection{Splice diagrams}\label{ssec:splice_diagrams}

We will recall the notions of maximal splice diagram and splice diagram. The maximal splice diagram \cite{NW05a} is a useful way to repackage the data of the plumbing graph $\Gamma$. It is directly related to the \emph{inverse matrix} of the plumbing matrix, which occurs in the formula for $\Zhat_{\sigma} (q)$. Even more importantly, it gives a good grasp of which entries of this matrix are important and which can be suppressed in some situations. This leads to the notion of the splice diagram. Splice diagrams are useful for constructing algebraic equations of singularities from topology; see Section~\ref{ssec:splice_quotient}.

\begin{Definition}
	The \emph{maximal splice diagram} of $Y$ is a graph on the same set of vertices as~$\Gamma$ with a weight $w_{ve}$ for each pair $(v, e)$ of an edge $e$ with its endpoint vertex $v$. The weight~$w_{ve}$ is given by the determinant of $-M^{ve}$, where $M^{ve}$ is the plumbing matrix of the connected component of $\Gamma{\setminus}e$ not containing $v$.

	The \emph{splice diagram} is formed from the maximal splice diagram by deleting the vertices of degree two and their weights, merging the neighboring edges, and finally removing the weights at leaves. In other words, we keep only the information of the weights around nodes.
\end{Definition}

\begin{Example}[\cite{NW05a}]\label{main_ex_H_shaped}
	Let $Y$ be a manifold given by the plumbing graph in Figure~\ref{Hshaped_plumbing}. It is a~homology sphere.
	Its maximal splice diagram and splice diagram are shown at Figure~\ref{Hshaped_splice_max}.
	\begin{figure}[ht]
 \centering
		\includegraphics[scale=1.0]{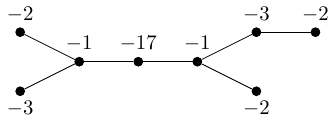}
		\caption{Plumbing graph of a homology sphere $Y$.}
		\label{Hshaped_plumbing}
	\end{figure}
	\begin{figure}[ht]
 \centering
		\includegraphics[scale=1.0]{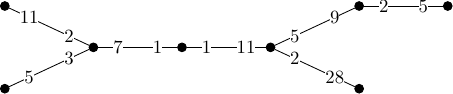}
		\hspace{0.3cm}
		\includegraphics[scale=1.0]{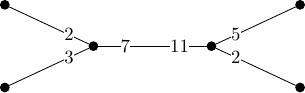}
		\caption{Maximal splice diagram and splice diagram of $Y$.}
		\label{Hshaped_splice_max}
	\end{figure}
\end{Example}

The plumbing graph and the maximal splice diagram are essentially equivalent as follows from the next theorem \cite{EN86}. For $v$, $w$ vertices of the maximal splice diagram, we define $N_{vw}$ as the product of all weights adjacent to the shortest path between $v$ and $w$, not lying on the path. For a matrix $A$, denote $\adj(A)$ the adjugate matrix of $A$, i.e., $\adj(A)=\det(A) A^{-1}$.
\begin{Theorem}[{\cite[Theorem~12.2]{NW05a}}]\label{theorem:adjugate} Let $M$ be a plumbing matrix of a negative definite manifold~$Y$ and let $N$ be as above. Then we have
$
		N = \adj(-M)$.
	Therefore, the maximal splice diagram of $Y$ determines the plumbing graph $\Gamma$ of $Y$.
\end{Theorem}

This theorem allows us to understand what information is contained in the splice diagram. It contains exactly the data of $\adj(-M)_{vw}$ for pairs $(v, w)$ of type (leaf, node), (node, leaf), (node, node), and finally (leaf, leaf) for two \emph{distinct} leaves.

The theorem is illustrated on Theorem~\ref{main_ex_H_shaped} in Figure~\ref{fig:path}.
	\begin{figure}[ht]
		\centering
		\begin{tikzpicture}
			\begin{scope}[every node/.style = {circle, fill, draw, inner sep = 1.5pt}]
				\node (1) at (-2.5, 0.7) {};
				\node (2) at (-2.5, -0.7) {};
				\node (3) at (-1, 0) {};
				\node (4) at (1, 0) {};
				\node (5) at (2.5, 0.7) {};
				\node (6) at (2.5, -0.7) {};
			\end{scope}

			\begin{scope}[every node/.style = {fill = white, inner sep = 0pt, scale = 0.9pt},
					every edge/.style = {draw}]
				\path[draw] (2) -- (3);
				\path[dashed] [-](1) edge node[near end] {$2$} (3);
				\path[draw] (3) -- (4);
				\path[dashed] [-](4) edge node[near start] {$5$} (5);
				\path[dashed] [-](4) edge node[near start] {$2$} (6);
			\end{scope}

			%vertex labels
			\node[xshift = -.3cm] at (2) {v};
			\node[yshift = .3cm] at (4) {w};

		\end{tikzpicture}
		\caption{The entry $\adj(-M)_{vw} = 2 \cdot 2 \cdot 5$ is given by the product of adjacent weights to the path connecting two vertices.}
		\label{fig:path}
	\end{figure}

\subsection{Universal abelian covers}

Recall our assumption that $Y$ is a rational homology sphere, so $H = H_1(Y, \ZZ)$ is finite with corresponding finite universal abelian cover $Y^{\rm ab}$. Universal abelian covers are related to splice diagrams by the following theorem, explaining their topological significance.

\begin{Theorem}[\cite{NW05a, Pe09}]\label{thm:ab_cover}
	Two manifolds\footnote{Satisfying our usual assumptions, see Section~\ref{ss:main_results}.} $Y_1$ and $Y_2$ have homeomorphic universal abelian covers if and only if they have the same splice diagram.
\end{Theorem}

It is often useful to consider $Y$ as the quotient of $Y^{\rm ab}$ by $H$ as quantities attached to $Y^{\rm ab}$ are often simpler, and then they can be subsequently equivariantly refined using the action of~$H$. In singularity theory, $Y^{\rm ab}$ is sometimes a link of a simpler singularity, e.g., it is of Brieskorn diagonal type for quasi-homogeneous singularities. One can think of $Z_0(q)$ as being attached to $Y^{\rm ab}$ and $\Zhat_\sigma(q)$ being the $H$-equivariant refinement. This is formalized in Theorem~\ref{thm:plumbing} and Corollary~\ref{cor:same_splice}.

Note that for $Y$ plumbed, the universal abelian cover $Y^{\rm ab}$ is plumbed as well. It does not have to be, however, a rational homology sphere, when $Y$ is, or some intermediate covering may lack this property. Unlike for the universal cover, taking the universal abelian cover is not an idempotent operation, i.e., as it may have nontrivial $H$ again and the construction can be repeated. For example, one has a tower of links of singularities
	\[
		\begin{tikzcd}
			S^3 \arrow[r, "{2:1}"] & A_1 \arrow[r, "{4:1}"] & D_4 \arrow[r, "{3:1}"] & E_6 \arrow[r, "{2:1}"] & E_7.
		\end{tikzcd}
	\]
It is an interesting question how the $\Zhat_\sigma(q)$ invariants behave with respect to taking covers. The only clue we have is the relation for lens spaces and their universal abelian cover, the 3-sphere~$S^3$, see Section~\ref{sec:lens}.

\subsection[Z_0 and splice diagrams]{$\boldsymbol{Z_0}$ and splice diagrams}

We will now prove that the series $Z_0(q)$ can be reconstructed from the splice diagram, up to a~prefactor, which is proportional to the Casson--Walker invariant $\lambda(Y)$ \cite{Wa92}. Let $H = H_1(Y, \ZZ)$ and $\lambda(Y)$ be the Casson--Walker invariant.
\begin{Theorem}\label{thm:plumbing}
	Let $Y$ be a negative definite plumbed manifold which is a rational homology sphere. Then the $q$-series
$
		q^{-6 \lambda(Y)} Z_{0}\bigl(q^{|H|}\bigr)
$
	only depends on the splice diagram of $Y$.
\end{Theorem}

\begin{proof}
	We will express the elements of the inverse matrix $M^{-1}$ using Theorem~\ref{theorem:adjugate}. Then we will find that those entries which are not expressible using splice diagram contribute only by an overall power of $q$. Here the integrand of $Z_0(q)$ plays an important role. Finally, we will relate this overall power to $\lambda(Y)$.
	
	The series $Z_0\bigl(q^{|H|}\bigr)$ can be written, up to prefactor $q^\Lambda$, as a sum of terms
	\begin{equation}\label{one_term}
		q^{-(\vl, \det(-M) M^{-1}\vl)/4} = q^{(\vl, \adj(-M)\vl)/4} = \prod_{v, w \in \Ver} q^{l_v l_w\adj(-M)_{vw}/4}
	\end{equation}
 for $\vl \in 2\mathbb{Z}^s+\vdelta$. Now if we look at the expansion of
	\begin{equation}
		\bigl( z_v - z_v^{-1}\bigr)^{2-\deg(v)}
	\end{equation}
	we see that a vector $\vl$ contributes to $Z_0\bigl(q^{|H|}\bigr)$ only if it has zero components for vertices of degree two. The only entries of $\adj(-M)$ which enter the formula are therefore $\adj(-M)_{vw}$ with $v$, $w$ leaves or nodes. By Theorem~\ref{theorem:adjugate}, all these entries are encoded in the splice diagram, except from $\adj(-M)_{vv}$ for a leaf $v$. In the product \eqref{one_term}, they contribute by
	\begin{equation}\label{one_term_power}
		\prod_{v \in \text{Leaves}} q^{\adj(-M)_{vv} l^2_v/4} = q^{\sum_{v \in \text{Leaves}} \adj(-M)_{vv}/4}
	\end{equation}
	because for leaves, we have $l_v^2 = (\pm 1)^2 = 1$, from the term $\bigl(z_v-z_v^{-1}\bigr)$. We see that these entries contribute to each monomial by the same overall power of $q$. Up to this power, $Z_0\bigl(q^{|H|}\bigr)$ depends only on the splice diagram.

	We are left to investigate the $q$-power
	\[
		\square := |H|(-3s-\Tr(M))/4 +\sum_{v \in \text{Leaves}} \adj(-M)_{vv}/4.
	\]
	A well-known formula for Casson--Walker invariant of negative definite plumbings $\lambda(Y)$ reads~\cite[p.\ 296]{NeNi02}
	\begin{equation}\label{eq:casson_formula}
		-\frac{24}{|H|}\lambda(Y) = \Tr(M) + 3s + \sum_{v \in \Ver} (2-\delta_v) M^{-1}_{vv},
	\end{equation}
 which can be rewritten using $|H|=\det(-M)$ as
 \begin{equation}\label{eq:casson_formula_adj}
		6\lambda(Y) = |H|(-3s-\Tr(M))/4 + \sum_{v \in \Ver} (2-\delta_v) \adj(-M)_{vv}/4,
	\end{equation}
	so
	\[
	\square = 6\lambda(Y) - \sum_{v \in \text{Nodes}} (2-\delta_v) \adj(-M)_{vv}/4.
	\]
	We see that the overall prefactor $\square$ equals $6\lambda(Y)$ plus terms that depend again only on the splice diagram.
\end{proof}

\begin{Remark}
	 The exponent $6\lambda(Y)/|H|$ in $Z_0(q)$ is the first coefficient $\lambda_1$ of the Ohtsuki series~\cite{Mu95,Oh96}. This is in agreement with the expectation that $Z_0(q)$ is the resummation of the Ohtsuki series \cite{GMP16}.

	The theorem formalizes the idea of the ``unimportance'' of the vertices of degree two in the plumbing graph $\Gamma$. It provides a computational tool for homology spheres where the unique~$\Zhat(q)$ and $Z_0(q)$ coincide, as the splice diagram is generally much smaller than the plumbing graph.
\end{Remark}

Using the relation of splice diagrams with universal abelian covers in Theorem~\ref{thm:ab_cover}, we obtain the following corollary, relating the series $Z_0(q)$ of two non-homeomorphic 3-manifolds.
\begin{Corollary}\label{cor:same_splice}
	If $Y_1$ and $Y_2$ have the same universal abelian cover, then
	\[
		q^{-6 \lambda(Y_1)} Z_0\bigl(Y_1, q^{|H_1(Y_1)|}\bigr) =
		q^{-6 \lambda(Y_2)} Z_0\bigl(Y_2, q^{|H_1(Y_2)|}\bigr).
	\]
\end{Corollary}
\begin{Example}
	Let us take the manifold from Theorem~\ref{main_ex_H_shaped} and call it $Y_1$. Another manifold $Y_2$ with the same splice diagram and $H_1 = \ZZ/17\ZZ$ is given in Figure~\ref{Hshaped_plumbing_17}.
	\begin{figure}
 \centering
		\includegraphics[scale=0.9]{graphics/main_example_H_shaped.pdf}
		\includegraphics[scale=0.9]{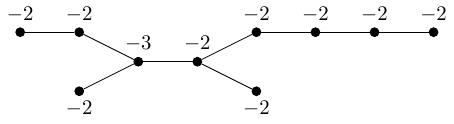}
		\caption{Two plumbed manifolds with the same splice diagram.}
 \label{Hshaped_plumbing_17}
	\end{figure}

	We have
	\begin{gather*}
		Z_{0}(Y_1, q) = \frac12 q^{7/2} \bigl(-1 + q - 2 q^2 + q^3 + q^5 + 3 q^9 + q^{10} - q^{14} - q^{16} -
		q^{17} +\cdots \bigr),\\
		Z_{0}\bigl(Y_2, q^{17}\bigr) = \frac12 q^{-53/2} \bigl(-1 + q - 2 q^2 + q^3 + q^5 + 3 q^9 + q^{10} - q^{14} - q^{16} -
		q^{17} + \cdots \bigr).
	\end{gather*}
	We see that they coincide up to an overall power. The Casson--Walker invariants are $\lambda(Y_1) = -4$ and $\lambda(Y_2) = -9,$ respectively.

	Note that the series $Z_0\bigl(Y_2,q^{17}\bigr)$ decomposes into 17 $q$-series $\Zhat_{\sigma} \bigl(Y_2, q^{17}\bigr)$. Since the $Z_0$ series for~$Y_1$ and $Y_2$ coincide up to the overall factor and scaling, we also have a decomposition of~$Z_0(Y_1,q)$. This decomposition would be far from obvious without the knowledge of $Y_2$ or, equivalently, the $\ZZ/17\ZZ$ action on $Y_1$.
\end{Example}

\begin{Remark}
	For rational homology spheres that are not Seifert-fibered, there is always a finite number of manifolds with a given splice diagram. Their number, however, grows rapidly \cite[Proposition~4.3]{Neu07}. When $Y$ is an integral homology sphere, then the splice diagram determines the maximal splice diagram and hence the plumbing graph \cite{EN86}.
\end{Remark}

\subsection{Splice and splice-quotient singularities}\label{ssec:splice_quotient}

In singularity theory, splice diagrams are used to define \emph{splice-quotient} singularities~\cite{NW05a}. This important class includes, e.g., weighted homogeneous, and rational singularities. Starting with a plumbed rational homology sphere $Y$ and its splice diagram, one first constructs the equations for the splice-type singularity\footnote{More precisely, an equisingular family of singularities.} $X^{\rm ab}$. These are equations in variables corresponding to the leaves of the splice diagram. They define an isolated complete intersection surface singularity. The splice-quotient $X$ is the quotient of $X^{\rm ab}$ by a natural action of $H_1(Y)$. Note that the construction relies on certain technical conditions on the plumbing graph.

The importance of these singularities stems from their topological nature. They provide a~good testing ground for comparing topological and analytical invariants. On the other hand, one should keep in mind that they are very special in some sense, even among the singularities with given link $Y$. A nice illustration of the above philosophy is provided by a theorem of N\'emethi.

\begin{Theorem}[{\cite[8.5.19--8.5.26]{Nem22}}]
	 A normal isolated surface singularity whose link is a rational homology sphere is a splice-quotient if and only if the analytical Poincar\'e series coincides with the topological Poincar\'e series.
\end{Theorem}

As we explained in Section~\ref{sec:plumbing}, $\Zhat_{\sigma} (q)$ can be thought of as $H$-equivariant decomposition of~$Z_0$. Hence we expect that they should be related to the splice-quotient $X$, whose link is the manifold~$Y$. The series $Z_0(q)$ should be related to the splice singularity, the universal abelian cover~$X^{\rm ab}$.
\begin{Question} Is there a complex analytic version of $\Zhat_{\sigma} (q)$ invariants such that it coincides with the usual $\Zhat_{\sigma} (q)$ in the case of quasi-homogeneous singularities, or more generally, splice-quotients?
\end{Question}
While we cannot provide the answer at the moment, we observe in Theorems~\ref{thm:plumbing} and~\ref{thm:seifert}, that the structural properties of $\Zhat_{\sigma} (q)$ are strikingly similar to the properties of invariants, such as the topological Poincar\'e series or Seiberg--Witten invariants, which have analytic analogs~\cite{Nem22}.

\section{Seifert manifolds}\label{sec:seifert}

We will now specialize our discussion to Seifert manifolds. For those, the plumbing graph can be chosen to be star-shaped, i.e., with at most one node. We assume that Seifert manifolds are fibered over $S^2$, are rational homology spheres and can be presented by a negative definite plumbing.

\begin{figure}
	\centering
	\includegraphics{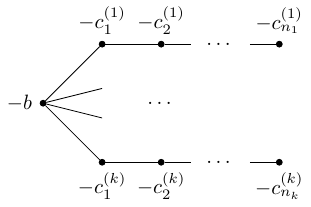}
	\caption{Plumbing graph of a Seifert manifold.}
	\label{fig:seifert}
\end{figure}

Let us fix some notation. We mostly follow \cite{Nem22} but use $b_i$ instead of $\omega_i$. In \cite{GM21}, the sign of $b$ and the roles of $a_i$ and $b_i$ are flipped. Seifert (reduced) data consist of an integer $b$ and tuples $(a_1, b_1), (a_2, b_2), \dots, (a_k, b_k)$ of integers such that $0<b_i<a_i$ and $\gcd(a_i, b_i) = 1$ for $i = 1,2,\dots,k$. Associated with it is a Seifert manifold $Y = M(b;(a_1, b_1), (a_2, b_2), \dots, (a_k, b_k))$ fibered over the sphere with $k$ singular fibers. We always assume that $k \geq 3$. Seifert manifolds with $k<3$ are lens spaces, and we can easily enlarge $k$ by adding tuples $(1, 0)$. The manifold $Y$ can be described by a star-shaped plumbing graph with the central node having the framing $-b$ and with $k$ legs, see Figure~\ref{fig:seifert}. The framing of vertices on $j$-th leg is given (starting at the node) by Hirzerbruch--Jung (HJ) continued fraction of $a_j/b_j$, namely the sequence \smash{$-c^{(j)}_1, -c^{(j)}_2, \dots, -c^{(j)}_{n_j}$} with \smash{$c^{(j)}_1 \geq 1$} and \smash{$c^{(j)}_i \geq 2$} for $i = 2, \dots, n_j$, such that
\[
	\frac{a_j}{b_j} = \big[c^{(j)}_1, c^{(j)}_2, \dots, c^{(j)}_{n_j}\big] := c^{(j)}_1 - \cfrac{1}{c^{(j)}_2-\cfrac{1}{\ddots - \cfrac{1}{c^{(j)}_{n_j}}}}.
\]

We label the leaves of the plumbing graph as $v_1, v_2, \dots, v_k$ and the node as $v_0$. We denote%
\[
A = \prod_{i = 1}^k a_i, \qquad A_i = \frac{A}{a_i}, \qquad A_{ij} = \frac{A}{a_ia_j}.
\]

The fundamental group $\pi_1(Y)$ has a presentation given by elements $g_1,g_2,\dots,g_k$ and $g_0$ satisfying
$[g_i, g_0] = 1$, $ g_i^{a_i} = g_0$, for $i = 1, 2, \dots, k$, $ \prod_{i = 1}^k g_i^{b_i} = g_0^b$.
The group $H = H_1(Y, \ZZ)$ is its abelianization, and we denote the generators by the same letters. Note that $g_i$ can be identified with the classes $[e_i] \in H = \ZZ^s/M\ZZ^s$ of the canonical basis vectors for the corresponding leaf $v_i$ and $g_0 = [e_0]$ is the class of the node.

Let
$
	e = -b + \sum_i \frac{b_i}{a_i}
$
be the Euler number of the Seifert fibration. Seifert manifold $Y$ over $S^2$ can be represented by a~negative definite plumbing if and only if $e<0$. Assuming this, $Y$ is a~rational homology sphere and the order of the group $H = H_1(Y, \ZZ)$ can be computed as
\[
	|H| = A |e| = \biggl|-A b + \sum_i A_i b_i\biggr|.
\]

\begin{figure}[ht]
	\centering
	\includegraphics{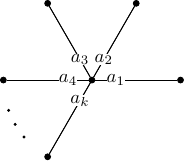}
	\caption{Splice diagram of Seifert manifold $M(b;(a_1, b_1), \dots, (a_k, b_k))$.}
	\label{fig:splice_seifert}
\end{figure}

\begin{Remark}
As we saw in Section~\ref{sec:splice}, some invariants of $Y$ can be read from the splice diagram of $Y$ while others need the full plumbing data. Our primary example are $Z_0(q)$ and $\Zhat_\sigma(q)$, respectively. In the Seifert case, the splice diagram contains exactly the integers $a_i$; see Figure~\ref{fig:splice_seifert}, so this distinction is about (in)dependence on $b$ and $b_i$.\footnote{And it is vacuous for integral homology spheres where $b$ and $b_i$ are uniquely determined by $a_i$.}

One has different examples of this phenomenon: While for fixed $a_1, a_2, \dots, a_k$, the order of $H$ can be arbitrarily large, the group $H/\langle g_0 \rangle$ is independent of $b_i$, with order $A/\lcm(a_1, a_2, \dots, a_k)$. Geometrically, $g_0$ is the generator corresponding to the fiber of the Seifert fibration, and $H/\langle g_0 \rangle$ is the fundamental group of the base orbifold $S^2$. Another interesting example is provided by the counts of $\operatorname{SL}_2(\C)$ vs $\operatorname{SU}(2)$ connections, see Section~\ref{sec:spectrum}.
\end{Remark}

\subsection{Reduction theorem for Seifert manifolds}

It is natural to ask if we can express the $\Zhat_\sigma(q)$ invariants of Seifert manifolds using Seifert data, rather than the plumbing graph (which is significantly larger due to continued fractions). This was done for Seifert homology spheres with three singular fibers (Brieskorn spheres) in \cite{GM21} and generalized slightly in \cite{Pa20a}, using false theta functions. Another approach to reduction of $\Zhat_{\sigma}$ appeared recently in \cite{Chung22}, for pairwise coprime $a_i$.

The main theorem of this section, Theorem~\ref{thm:seifert}, generalizes these results to any number of fibers, removes the conditions on $H$ and, perhaps more importantly, makes the role of Seifert data more transparent. It also emphasizes the role of the series $Z_0(q)$, treating $\Zhat_{\sigma} (q)$ as its $H$-equivariant refinement, using the anticanonical \spc{} structure $\ac$ from Section~\ref{spc_structures}. This ``reduction theorem'' gives us a formula for $Z_0(q)$ and $\Zhat_\sigma(q)$ using single variable rational functions. Recall the ``Laplace transform'' $\LC_A$ defined as
\begin{equation}\label{laplace}
	\LC_A(t^n) = q^{n^2/4A}.
\end{equation}

\begin{Theorem}[reduction theorem]\label{thm:seifert}
	Let $Y = M(b;(a_1, b_1), (a_2, b_2), \dots, (a_k, b_k))$ be a Seifert manifold over $S^2$ $(k \geq 3)$ with $e<0$. Then
$
		Z_0\bigl(q^{|H|}\bigr) = q^\Lambda \LC_A(\se f_0(t))$,
	where
	\begin{equation}\label{seifert}
		f_0(t) = \frac{\bigl(t^{A_1}-t^{-A_1}\bigr)\bigl(t^{A_2}-t^{-A_2}\bigr)\cdots\bigl(t^{A_k}-t^{-A_k}\bigr)}{\bigl(t^A-t^{-A}\bigr)^{k-2}},
	\end{equation}
	and $\Lambda = \Lambda(Y)$ is a rational number described in the proof.
	Moreover, we have
	\[
		\Zhat_{h \ac}(q) = q^{\Lambda/|H|} \LC_{A|H|}(\se f_{h \ac}(t)),
	\]
	where
	\begin{equation}
		\sum_{h \in H} f_{h \ac}(t)h = t^{(k-2)A - \sum_i A_i}\frac{\bigl(g_1 t^{2A_1}-1\bigr)\cdots\bigl(g_k t^{2A_k}-1\bigr)}{\bigl(g_0 t^{2A}-1\bigr)^{k-2}}.
	\end{equation}
\end{Theorem}

\begin{proof}
	We will start with the series $Z_0(q)$ for simplicity. The idea is to reduce the multivariate integrand of $Z_0$ to univariate rational function $f_0(t)$, and then study the effect of this operation on the exponents given by the theta function. The function $f_0(t)$ is obtained from the integrand in $Z_0(q)$ by the following substitution
	\begin{gather}\label{subst}
			z_i \ra t^{A_i} \quad\mbox{for the leaves},\qquad
			z_0 \ra t^A \quad \mbox{for the node}.
	\end{gather}
	When computing $Z_0\bigl(q^{|H|}\bigr)$ from the definition, we substitute each monomial $\vz ^{\vl}$ with \smash{$q^{(\vl, \adj(-M)\vl)\!/\!4}$}.
 We have the following relations from Theorem~\ref{theorem:adjugate} (see Figure~\ref{fig:splice_seifert})
 \begin{equation}\label{eq:adjugate}
 \adj(-M)_{00} = A, \qquad \adj(-M)_{0i} = A_i, \qquad \adj(-M)_{ij} = A_{ij}\qquad \text{for}\quad i \neq j.
 \end{equation}
 We can use \eqref{eq:adjugate} to expand the exponent (ignoring factors of 4 in what follows)
	\begin{equation}\label{first}
		\bigl(\vl, \adj(-M)\vl\bigr) = l_0^2 A+2\sum_{i\neq 0} l_0l_i A_{i} + \sum_{\substack{i \neq j\\ i,j \neq 0}} l_i l_j A_{ij}+\sum_{i\neq 0} l_i^2 \adj(-M)_{ii}.
	\end{equation}
	On the other hand, if we first perform the substitution \eqref{subst} and then Laplace transform~$\LC_A$~\eqref{laplace}, we get
	\begin{align}
\biggl(l_0 A + \sum_{i\neq 0}l_i A_i \biggr)^2/A & = l_0^2 A + 2 \sum_{i\neq 0} l_0 l_i A_i + \sum_{i, j\neq 0} l_i l_j A_{i}A_{j}/A \nonumber \\
& = l_0^2 A + 2 \sum_{i\neq 0} l_0 l_i A_i + \sum_{i, j\neq 0} l_i l_j A_{ij}.\label{second}
 \end{align}
We see that the two expressions \eqref{first} and \eqref{second} differ only in terms corresponding to ``a leaf with itself'', $l_i^2 A_{ii}$ vs. $l_i^2 \adj(-M)_{ii}$. Since $l_i^2 = (\pm 1)^2 = 1$ for each vector $\vl \in \ZZ^s$ that contributes a nonzero term to $Z_0(q)$, these numbers give the same contribution to every term of the sum, so they contribute by an overall power of $q$ in the series $Z_0(q)$,
 \begin{align}
		4\Lambda
		 & = |H| (-3s - \Tr(M)) + \sum_{i\neq 0}\adj(-M)_{ii} - \frac{1}{A}\sum_{i\neq 0} A_i^2 \nonumber\\
		 & = -|H|(3s+\Tr(M))+\bigg(\sum_{i}(2-\delta_i)\adj(-M)_{ii}\bigg) - (2-\delta_0) \adj(-M)_{00} - \frac{1}{A}\sum_{i\neq 0} A_i^2\nonumber\\
 & = -|H|\bigg(3s+\Tr(M)+\sum_{i}(2-\delta_i)(M^{-1})_{ii}\bigg) - (2-\delta_0) A - \frac{1}{A}\sum_{i\neq 0} A_i^2\nonumber\\
		 & = 24 \lambda(Y) - A\bigg(2-\delta_0 + \sum_{i\neq 0} \frac{1}{a_i^2}\bigg).\label{eq:Lambda}
 \end{align}
	In the last equation, we used the formula \eqref{eq:casson_formula} for Casson--Walker invariant $\lambda$ as in the proof of Theorem~\ref{thm:plumbing}.
	To obtain the formula for $\Zhat_\sigma(q)$, we follow the same computation of the exponents (without scaling $q$ by $|H|$), applied on the equivariant expression \eqref{eq:zhat_expansion_can}.
\end{proof}

We now give some remarks on interpreting this result and its relation to previous works. Seifert manifolds have often been used as examples for exploring interesting phenomena for~$\Zhat(q)$, e.g., \cite{CCFGH18,CCFFGHP22, CFS20, MT21, Sug22}. Our theorem gives an easy-to-use formula that can be used to check many conjectures, such as connections with logarithmic vertex algebras in Section~\ref{sec:spectrum}.

The reduction theorem can be considered an analogue of the reduction procedure of N\'emethi \cite[p.~364]{Nem22}, which was developed in greater generality. It is natural to proceed in this direction for more general graphs and try to obtain more general reduction theorems for $\Zhat_{\sigma} (q)$ invariants. A clue toward a more general reduction was already presented in \cite{BMM20} for $H$-shaped graphs with exactly 6 nodes. Splice diagrams provide a natural framework for generalizing these computations for any plumbing.

The function $f_0(t)$ is a ``symmetrized and inverted version'' of the univariate Poincar\'e series of the splice type singularity defined by the splice diagram. For Seifert manifolds, those are complete intersections of Brieskorn type \cite{Neu83}. Following \cite{AM22,GMP16}, we can also identify the function~$f_0(t)$ with the Borel transform of the perturbative expansion in complex Chern--Simons theory; its singularity structure is the central element of the resurgent analysis. Indeed, as explained in the Introduction, $Z_0$ is a very natural object from the viewpoint of complex Chern--Simons theory and resurgent analysis. Unfortunately, it generally lacks integrality and many other important properties, e.g., it is not expected to admit a categorification, whereas $\Zhat_{\sigma}(q)$ enjoy these properties. Our discussion here seems to suggest that for the latter, the role of the Borel plane is played by the Poincar\'e series of the corresponding splice-type singularity. It would be interesting to explore this further.

A useful corollary of Theorem~\ref{thm:seifert} explains frequent vanishing of the $\Zhat_{\sigma}(q)$ series for Seifert manifolds.

\begin{Corollary}\label{cor:vanishing}
	Let $Y$ be a Seifert manifold with $k$ singular fibers as above, for which the element $g_0 \in H_1(Y,\ZZ)$ is trivial, i.e., $|e| \lcm(a_1,\dots,a_k) =1$. Then there are at most $2^k$ nonzero~$\Zhat_{\sigma}$ invariants.
\end{Corollary}

\begin{proof}
	As $g_0$ is trivial, the only elements of $H$ which appear in the expansion of
	\begin{equation}\label{seifert_equiv}
		\sum_h f_{h \cdot \ac}(t)h = t^{(\dots)}\frac{ \bigl(g_1 t^{2A_1}-1\bigr)\cdots\bigl(g_k t^{2A_k}-1\bigr)}{\bigl(g_0t^{2A}-1\bigr)^{k-2}}
	\end{equation}
	are the products of the generators in the numerator. Hence each such element $h$ is of the form~${
		h = g_1^{u_1} g_2^{u_2} \cdots g_k^{u_k}}$,
	where $u_i$ is $0$ or $1$. There are at most $2^k$ such products, so at most $2^k$ elements $h$ with a nonzero coefficient $f_{h \cdot \ac}(t)$.
\end{proof}

\begin{Example}
	For any positive integer $n$, Seifert manifold $M(2;(n, n-1), (n, n-1), (n, 1))$ has $|H| = n^3/n = n^2$ and $|g_0| = n/n = 1$ so there are only eight nonzero $q$-series $\Zhat_{\sigma} (q)$.
	On the other hand, consider $M(1;(3, 1), (4, 1), (5, 1))$.\footnote{The link of $S_{12}$ singularity $x^3y+y^2z+xz^2 = 0$.} It has $|H| = 13$ and all thirteen $\Zhat_{\sigma} (q)$ are nonzero, as $g_0$ generates $H$.
\end{Example}

\begin{Example}\label{ex:E6}
	Let us illustrate in detail how the reduction theorem works on the Seifert manifold $Y = M(2;(2, 1), (3, 2), (3, 2))$, the link of $E_6$ singularity $x^2+y^3+z^4 = 0$. A plumbing graph and splice diagram for $Y$ are given in Figure~\ref{fig:E6_plumbing}.
	\begin{figure}
 \centering
		\includegraphics{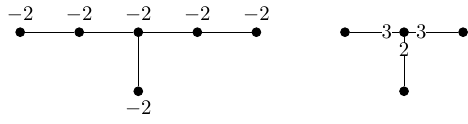}
		\caption{Plumbing graph and splice diagram of Seifert manifold $Y = M(2;(2, 1), (3, 2), (3, 2))$.}
		\label{fig:E6_plumbing}
	\end{figure}

	$Z_0(q)$: Consider the function $f_0$
	\begin{equation}
		f_0(t) = \frac{\bigl(t^{9}-t^{-9}\bigr)\bigl(t^{6}-t^{-6}\bigr)\bigl(t^{6}-t^{-6}\bigr)}{t^{18}-t^{-18}} = \frac{\bigl(t^{9}-t^{-9}\bigr)\bigl(t^{12}+t^{-12}-2\bigr)}{t^{18}-t^{-18}}.
	\end{equation}
	If we expand $f_0(t)$ as a power series in $t^{\pm}$ using the symmetric expansion, we obtain
	\[
	\frac{1}{2}\bigl(\dots +t^{-33} + 2 t^{-27} + t^{-21} - t^{-15} - 2 t^{-9} + t^{-3} + t^{3} - 2 t^{9} - t^{15} + t^{21} + 2 t^{27} + t^{33} + \cdots\bigr).
	\]
	Laplace transform $\LC_{18}$ substitutes each monomial $t^n$ with $q^{n^2/{4 \cdot 18}}$
	\[
		q^{1/8} \bigl(1 - 2 q - q^3 + q^6 + 2 q^{10} + q^{15} - q^{21} - 2 q^{28} - q^{36} +
	 q^{45} + 2 q^{55} +\cdots \bigr) = q^{25/8} Z_0\bigl(q^3\bigr).
	\]
	The exponent $25/8$ is related to the Casson--Walker invariant $\lambda(Y) = -11/12$ as in the proof
	\[
		 - \frac{25}{8} = -\Lambda = 6 \lambda(Y) + \frac{18}{4} -\frac{6^2 + 6^2 + 9^2}{4 \cdot 18}.
	\]
	$\Zhat_\sigma(q)$: We have $e = -1/6$ and $|H| = |e| a_1 a_2 a_3 = 3$. The group $H$ is generated by $g_1$, $g_2$, $g_3$, $g_0$ satisfying $g_1^2 = g_0$, $g_2^3 = g_0$, $g_3^3 = g_0$, $g_1g_2^2g_3^2 = g_0^{2}$. It is easy to see that $g_0 = g_1 = 1$, and put $g := g_2 = g_3^2$ and $\sigma := \ac$. The rational functions in Theorem~\ref{thm:seifert} read
	%The chern class of $\can$ is $g_1 g_2 g_3 g_0^{-1}=g \cdot g^2 = 1$ so it is spin.
	\begin{align*}
		\sum_h f_{h \cdot \sigma}(t) h &= t^{18-9-6-6} \frac{\bigl(g_1 t^{18}-1\bigr)\bigl(g_2 t^{12}- 1\bigr)\bigl(g_3 t^{12}-1\bigr)}{g_0 t^{36}-1} = t^{-3}\frac{\bigl(g t^{12}- 1\bigr)\bigl(g^2 t^{12}-1\bigr)}{1+t^{18}}\\
							&= \frac{t^{21}+t^{-3}}{1+t^{18}} + \bigl(g+g^2\bigr)\frac{-t^9}{1+t^{18}}.
	\end{align*}
	These functions are expanded as
	\begin{gather*}
		\se f_{\sigma}(t) = \frac12 \bigl(\dots + t^{-21}-t^{-15}+t^{-3}+t^3-t^{-15}+t^{21} + \cdots \bigr),\\
		\se f_{g\sigma}(t) = \se f_{g^2\sigma}= \frac12 \bigl(\dots-t^{-45}+t^{-27}-t^{-9}-t^{9}+t^{27}-t^{45} + \cdots \bigr).
	\end{gather*}
	After performing the Laplace transform $\LC_{A |H|}\colon t^n \to q^{n^2/4\cdot 18\cdot 3}$ and multiplying by $q^{\Lambda/3}=q^{25/24}$ we obtain the $q$-series:
	\begin{gather*}
		\Zhat_{\sigma}(q) = q^{-1} - 1 + q + q^4 - q^6 - q^{11} + q^{14} + q^{21} - q^{25} - q^{34} + q^{39} + q^{50} - q^{56} - \cdots, \\
		\Zhat_{g\sigma}(q) = \Zhat_{g^2\sigma}(q) = q^{-2/3}\bigl(-1 + q^3 - q^9 + q^{18} - q^{30} + q^{45} - q^{63} + q^{84} - q^{108} + \cdots \bigr).
	\end{gather*}
\end{Example}

\subsection{Lens spaces}\label{sec:lens}

In the previous subsection, we studied Seifert manifolds with 3 or more singular fiberes, but the method of the proof of Theorem~\ref{thm:seifert} also lets us to compute $\Zhat_\sigma(q)$ of lens spaces $L(p, r)$. This has been done systematically only in the case of $r = 1$ in \cite{GPP20}. The lens space $Y = L(p,r)$ with~${p>r>0}$ can be expressed as a plumbed manifold given by a path
\begin{center}
\includegraphics{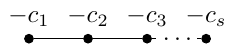}
\end{center}
where the coefficients are given by the continued fraction $p/r = [c_1,c_2,\dots,c_s]$.
The rational function in the integrand in \eqref{eq:plumbing1} is simply
$
 \bigl(z_1-z_1^{-1}\bigr)\bigl(z_s-z_s^{-1}\bigr)$,
with two variables corresponding to the leaves. The relevant entries of the adjugate of the plumbing matrix are
$\adj(-M)_{11}=r$, $ \adj(-M)_{ss}=r^*$, $ \adj(-M)_{1s}=1$,
where $0<r^*<p$ satisfies $rr^* \equiv 1 \pmod{p}$.

Following the idea of the proof of Theorem~\ref{thm:seifert}, we can set $z_1=t$ and $z_s=t^r$ and express~$Z_0(q^p)$ as a Laplace transform in a single variable $t$. The overall $q$-power is fixed by comparing analogues of~\eqref{first} and~\eqref{second},
\begin{align*}
	Z_0(Y; q^p) & =q^{p\cdot\frac{-3s-\Tr(M)}{4}}q^{\adj(-M)_{11}+\adj(-M)_{ss}+2} q^{-\frac{r}{4}-\frac{1}{4r}-2} \LC_r\bigl(\bigl(t-t^{-1}\bigr)\bigl(t^r-t^{-r}\bigr)\bigr)\\
	 & = 2q^{6\lambda(Y)}q^{-\frac{r}{4}-\frac{1}{4r}} \LC_r \bigl(t^{r+1}-t^{r-1}\bigr)= 2q^{6\lambda(Y)} \bigl(q^{1/2}-q^{-1/2}\bigr) = q^{6\lambda(Y)} Z_0\bigl(S^3, q\bigr).
\end{align*}

In the second equation, we used again the formula~\eqref{eq:casson_formula_adj} for the Casson--Walker invariant~$\lambda(Y)$, which can also be expressed as $p \cdot \ds(r,p)/2$ where $\ds(r,p)$ is the Dedekind sum \cite[Example~9.3.6]{Nem22}. The independent factor $\bigl(q^{1/2}-q^{-1/2}\bigr)$ is the $Z_0(q)$ of the 3-sphere $S^3$. As $S^3$ is the universal (abelian) cover of the lens space, we get a direct relation between the invariants of the base space and the cover.

For the equivariant version, we have
\begin{align}
	\sum_{h \in H} \Zhat_{h \ac}(q)h & = q^{6\lambda(Y)}q^{-\frac{r}{4}-\frac{1}{4r}} \LC_b \bigl(t^{-1-r}\bigl(g_1 t^2-1\bigr)\bigl(g_s t^{2r}-1\bigr)\bigr)\nonumber\\
	 & = q^{6\lambda(Y)/p} \bigl((g_1g_s+1) q^{1/2p} - (g_1+g_s) q^{-1/2p}\bigr)\nonumber\\
				& = q^{6\lambda(Y)/p} \bigl(\bigl(g_1^{r+1}+1\bigr) q^{1/2p} - \bigl(g_1^{r}+g_1\bigr) q^{-1/2p}\bigr),\label{lens_equiv}
\end{align}
where $g_1$ and $g_s=g_1^r$ are the generators of $H \cong \Z_p$ corresponding to the leaves.

As a corollary, the collection \smash{$\big\{\Zhat_{\sigma}(q)\big\}_{\sigma \in \spc(Y)}$} distinguishes lens spaces. Indeed, the expression \ref{lens_equiv} specializes for $q=1$ to the Reidemeister--Turaev torsion (as conjectured in \cite{CGPS20}), which is known to have this property. It is interesting to notice that one needs to consider the whole package of $\Zhat_\sigma(q)$ invariants --- the set of exponents of $q$ itself does not distinguish lens spaces as there are distinct lens spaces with the same Casson--Walker invariant.

\section[Z\_0, irreducible connections, spectrum]{$\boldsymbol{Z_0}$, irreducible connections, spectrum}\label{sec:spectrum}

In this last section, we explore the interplay between four different objects associated with a~Seifert manifold $Y$: the $q$-series $Z_0$, irreducible flat $\operatorname{SL}_2(\C)$-connections, the spectrum of the associated splice-type singularity, and the related vertex algebras. In Section~\ref{sec:splice}, we examined which invariants of the 3-manifold depend on the full data of the plumbing graph, or rather just on the splice diagram. With the exception of vertex algebras, all objects listed above are of the latter type.\footnote{In particular, all of them become essentially trivial for lens spaces.} It is an interesting question whether the discussion in this section generalizes to more complicated plumbed manifolds.

This connection also hints at a deeper relation between geometric structures, leading to a~conjectural notion of a nonabelian noncommutative Hodge structure. Seifert manifolds form a~rather special class, and their properties are closely related to the base orbifold curve. Some of the relations described may be of algebraic origin rather than topological.

\subsection{Irreducible flat connections on Seifert manifolds and spectrum}

In this subsection we study the relation of flat connections on Seifert manifolds and the spectrum of the corresponding splice-type singularities. We extend the relation between the $\operatorname{SL}_2(\C)$ Casson invariant and Milnor number of Brieskorn spheres \cite[equation~(2.9)]{BC06} in two directions --- to Seifert rational homology spheres with 3 singular fibers, and to Seifert homology spheres with an arbitrary number of singular fibers. The argument is elementary, based on previously known results, but putting them together requires some care.

We first recall the case of Brieskorn spheres. For pairwise coprime integers $a_1,a_2,a_3 \geq 2$, consider the equation
\begin{equation}\label{brieskorn}
 X_t\colon\ x_1^{a_1}+x_2^{a_2}+x_3^{a_3}=t.
\end{equation}
For $t=0$, the hypersurface $X=X_0$ defines an isolated singularity at the origin whose link is the Brieskorn sphere $Y = \Sigma(a_1,a_2,a_3)$. The family $X_t$ is a natural one-parameter smoothing of $X$ with the corresponding Milnor fiber $F$ and the monodromy operator $h$ on $H^2(F,\C)$. The eigenvalues of $h$ are roots of unity. Moreover, there is a canonical mixed Hodge structure on the cohomology of $F$ which allows us to lift the eigenvalues of $h$ to rational numbers via the exponential map $x \to \exp(2\pi {\rm i} x)$. They are called \emph{spectral numbers} and their set is the \emph{spectrum} of the singularity $X$.

The spectrum of \eqref{brieskorn} has the following simple description. Let $D_3$ be the collection of tuples~${(d_1,d_2,d_3)}$ satisfying $0<d_i<a_i \text{ for } i = 1, 2, 3$. Then the spectral numbers are
\begin{equation}\label{spec_three}
	\frac{d_1}{a_1} + \frac{d_2}{a_2} + \frac{d_3}{a_3}, \qquad (d_1,d_2,d_3) \in D_3.
\end{equation}
Their number equals the total rank of $H^2(F,\C)$, the Milnor number $\mu(X)$,
$
\mu(X) = |D_3| = (a_1-1)(a_2-1)(a_3-1)$.

The $\operatorname{SL}_2(\C)$ Casson invariant \smash{$\lambda^C_{\operatorname{SL}_2(\C)}(Y)$} was defined in \cite{Cu01}. For Brieskorn spheres, Boden and Curtis proved that it is equal to the number of characters of irreducible flat $\operatorname{SL}_2(\C)$-connections \cite{BC06}. Those can be labeled by the orbits $D_3/(\Z_2)^2$, where the free $(\Z_2)^2$-action is given by the cyclic permutations of
\begin{equation}\label{eq:action}
d_1 \to a_1 - d_1, \qquad d_2 \to a_2 - d_2, \qquad d_3 \to d_3.
\end{equation}
This gives the following.

\begin{Theorem}[{\cite[Theorem~2.4]{BC06}}]\label{CassonMilnor}
With the notation above,
\smash{$\lambda^C_{\operatorname{SL}_2(\C)}(Y) = \frac{\mu(X)}{4}$}.
\end{Theorem}

\subsubsection{Rational homology spheres}

Let $Y = M(b;(a_1, b_1), \dots, (a_3,b_3))$ be a Seifert rational homology sphere. To obtain a generalization of Theorem~\ref{CassonMilnor} for this case, we need to modify both sides of the equation. On the singularity side, the equation \eqref{brieskorn} is the splice-type equation constructed from the splice diagram of $Y$, or equivalently, the universal abelian cover $X^{\rm ab}$ of a quasihomogeneous singularity~$X$ with the link $Y$ \cite{Neu83}. Moreover, the $(\Z_2)^2$-action \eqref{eq:action} on the spectrum of $X^{\rm ab}$ is no longer free. We denote by $\bigl\lvert\spec\bigl(X^{\rm ab}\bigr)/(\Z_2)^2\bigr\rvert$ the number of orbits of the $(\Z_2)^2$-action \eqref{eq:action} on \eqref{spec_three}.

On the topology side, we need to consider characters of all non-abelian flat $\operatorname{SL}_2(\C)$ connections, as opposed to the irreducible ones. Their number was computed by Cui, Qiu and Wang~\cite{CuQiWa21}. We denote it by $\lambda^{\rm nab}_{\operatorname{SL}_2(\C)}(Y)$ in an analogy with the $\operatorname{SL}_2(\C)$ Casson invariant.

\begin{Proposition}\label{prop:rat_sph_spec}
With the notation above,
\smash{$\lambda^{\rm nab}_{\operatorname{SL}_2(\C)}(Y) =\bigl\lvert\frac{\spec(X^{\rm ab})}{(\Z_2)^2}\bigr\rvert$}.
\end{Proposition}

\begin{proof}
 The formula for \smash{$\lambda^{\rm nab}_{\operatorname{SL}_2(\C)}(Y)$} in \cite[p.~20]{CuQiWa21} reads
 \begin{equation}\label{eq:conn_count_rat}
\lambda^{\rm nab}_{\operatorname{SL}_2(\C)}(Y) = \floor*{\frac{a_1}{2}}\floor*{\frac{a_2}{2}}\floor*{\frac{a_3}{2}}+\Big\lfloor\frac{a_1-1}{2}\Big\rfloor\Big\lfloor\frac{a_2-1}{2}\Big\rfloor
 \Big\lfloor\frac{a_3-1}{2}\Big\rfloor.
 \end{equation}
 The orbits of the action \eqref{eq:action} can be computed using Burnside's lemma
 \begin{equation}\label{eq:burside}
 \big\lvert D_3/(\Z_2)^2\big\rvert = \frac{1}{\big\lvert(\Z_2)^2\big\rvert} \sum_{g \in (\Z_2)^2} D_3^g.
 \end{equation}
 Here $D_3^g$ denotes the set of tuples in $D_3$ fixed by an element $g \in (\Z_2)^2$.
Denote
 \[
 \bar{a} =
 \begin{cases}
 1 & \text{if } a \text{ is even},\\
 0 & \text{if } a \text{ is odd}.\\
 \end{cases}
 \]
Then from \eqref{eq:burside} we have
\begin{align}
\big\lvert D_3/(\Z_2)^2\big\rvert
={}& \frac14 (\bar{a}_1\bar{a}_2(a_3-1)+\bar{a}_2\bar{a}_3(a_1-1)\nonumber\\
&+\bar{a}_3\bar{a}_1(a_2-1)+(a_1-1)(a_2-1)(a_3-1)).\label{eq:count_orbits}
\end{align}
We can write \eqref{eq:conn_count_rat} as
\[
 \left(\frac{a_1-1}{2}+\frac{\bar{a}_1}{2}\right)\left(\frac{a_2-1}{2}+\frac{\bar{a}_2}{2}\right)\left(\frac{a_3-1}{2}+\frac{\bar{a}_3}{2}\right)+\left(\frac{a_1-1}{2}-\frac{\bar{a}_1}{2}\right)\cdots.
\]
We see that the terms containing an even number of factors of the type $\frac{a_i-1}{2}$ cancel out and we obtain the expression \eqref{eq:count_orbits}.
\end{proof}

\subsubsection{Seifert homology spheres}

Let $Y = M(b;(a_1, b_1), \dots, (a_k, b_k))$ be a Seifert homology sphere with $k \geq 3$ singular fibers. The corresponding splice-type singularity is the Brieskorn--Hamm complete intersection \cite{Neu83}. It is given by a system of $k-2$ equations in $k$ variables of the following form: Let $\alpha_{ij}$ be a~complex matrix of size $(k-2) \times k$ with all maximal subdeterminants being nonzero. Consider the functions~$f_i(x)\colon \CC^k \rightarrow \CC$, $i = 1, \dots, k-2$, given by
\begin{equation}\label{icis}
	f_i(x) = \sum_{j=1}^k \alpha_{ij} x_j^{a_j}.
\end{equation}
The system
$
	f_1(x) = f_2(x) = \dots = f_{k-2}(x) = 0
$
defines a complex surface $X$ with an isolated complete intersection singularity at the origin, whose link is $Y$.

Let us denote by \smash{$\mathcal{M}^*(Y, \operatorname{SL}_2(\C))$} the moduli space of irreducible flat $\operatorname{SL}_2(\C)$ connections on $Y$ modulo conjugation. This moduli space can be described as a certain moduli space of parabolic Higgs bundles over the (orbifold) base $\CC\PP^1$ \cite{BoYo96}. It follows that \smash{$\mathcal{M}^*(Y, \operatorname{SL}_2(\C))$} is smooth, and has components of complex dimensions $0,2,4, \dots, 2k-6$.

The $\operatorname{SL}_2(\C)$ Casson invariant \smash{$\lambda^C_{\operatorname{SL}_2(\C)}(Y)$} counts the zero-dimensional components of $\mathcal{M}^*(Y, \allowbreak \operatorname{SL}_2(\C))$.
An alternative definition, denoted by \smash{$\lambda^P_{\operatorname{SL}_2(\C)}(Y)$}, is due to Abouzaid and Manolescu~\cite{AbMa20}. For a Seifert homology sphere $Y$, it coincides with the Euler characteristic of $\mathcal{M}^*(Y, \operatorname{SL}_2(\C))$. We will demonstrate that this version of $\operatorname{SL}_2(\C)$ Casson invariant is related to the Milnor number of the complete intersection $X$.

\begin{Theorem}\label{thm:seifert_casson}
	Let $Y$ be a Seifert homology sphere and $X$ the corresponding Brieskorn--Hamm complete intersection singularity with Milnor number $\mu(X)$. Then
\smash{$\lambda^{P}_{\operatorname{SL}_2(\C)}(Y) = \frac{\mu(X)}{4}$}.	
\end{Theorem}

\begin{proof}

The $\operatorname{SL}_2(\C)$-representations of $\pi_1(Y)$ were studied in \cite[proof of Theorem 2.7]{BC06}. Here we present an equivalent description, similar to the one given in \cite{AM22}.\footnote{Note that some components of lower dimensions were omitted there due to a too restrictive fundamental domain for the $(\ZZ_2)^{n-1}$ action below.}
For any $n$ satisfying~${3 \leq n \leq k}$, the set of components of $\mathcal{M}^*(Y, \operatorname{SL}_2(\C))$ of dimension $2n-6$ can be identified with the set of orbits
$
D_n/(\Z_2)^{n-1}$.
Here $D_n$ consists of tuples $(d_1, d_2, \dots, d_k)$ satisfying $0<d_i<a_i$ with exactly~$n$ numbers $d_i$ being nonzero. The action of $(\Z_2)^{n-1}$ takes $d_i \rightarrow a_i - d_i$ simultaneously for an even number of nonzero $d_i$, generalizing \eqref{eq:action}. The numbers $d_i$ can be thought of as rotational numbers of the corresponding representation with the choices of a representative in the $(\Z_2)^{n-1}$ orbit corresponding to taking different presentations of the fundamental group.

We obtain that the number of components of the dimension $2n-6$ is
\begin{equation}\label{count_rep}
	\frac{e_n(a_1-1, a_2-1, \dots, a_k-1)}{2^{n-1}},
\end{equation}
where $e_n$ is the elementary symmetric polynomial
\[
	e_n(x_1, x_2, \dots, x_k) = \sum_{1\leq j_1 < j_2 <\dots <j_n \leq k} x_{j_1}x_{j_2}\cdots x_{j_n}.
\]

Boden and Yokogawa \cite{BoYo96} computed the Poincar\'e polynomial of the components of the moduli space $\mathcal{M}^*(Y, \operatorname{SL}_2(\C))$. In particular, for any component $C \subset \mathcal{M}^*(Y, \operatorname{SL}_2(\C))$ of dimension $2n-6$, its Euler characteristic is
\begin{equation}\label{euler_char}
	\chi(C) = (n-1)(n-2)2^{n-4} = \binom{n-1}{2}2^{n-3}.
\end{equation}

From \eqref{count_rep} and \eqref{euler_char}, we obtain that the Euler characteristic of $\mathcal{M}^*(Y, \operatorname{SL}_2(\C))$ is
\begin{align}
\chi(\mathcal{M}^*(Y, \operatorname{SL}_2(\C))) &=
 \sum_{3 \leq n \leq k } \frac{e_n(a_1-1, \dots, a_k-1)}{2^{n-1}} \binom{n-1}{2}2^{n-3}\nonumber \\
 &= \frac14 \sum_{3 \leq n \leq k } e_n(a_1-1, \dots, a_k-1) \binom{n-1}{2}.\label{eq:euler_char_moduli}
\end{align}

 We now describe the monodromy eigenvalues (and hence the Milnor number), as computed by Hamm \cite{Ham72}. Let $X^*$ be the surface given by the equations \eqref{icis} except the last one ($f_{k-2}=0$). The function $f_{k-2}$ is regular on $X^*$ except at $0$, and we can define the monodromy operator $h$ of the Milnor fiber of $f_{k-2}$. We have the following.

\begin{Lemma}[{\cite[Lemma 4.2]{Ham72}}]
	The characteristic polynomial $p(t)$ of the monodromy $h$ is given~by
	\[
		p(t) = \prod_{\substack{3 \leq n \leq k, \\ 1 \leq j_1 < \dots < j_n \leq k }} \delta_{a_{j_1}, a_{j_2}, \dots, a_{j_n}}(t)^{\binom{n-1}{2}},
	\]
	where
	\[
		\delta_{c_1, c_2, \dots, c_n}(t) = \prod_{\substack{d_1, \dots, d_n, \\ 1 \leq d_i < c_i, \\ i = 1, \dots, n}}
		\bigl({\rm e}^{2\pi {\rm i} (d_1/c_1+d_2/c_2+\cdots+d_n/c_n)}-t\bigr).
	\]
\end{Lemma}

We see that the eigenvalues of $h$ are labeled by tuples $(d_1, d_2, \dots, d_n)$ for $n$ from 3 to $k$ in the same way as are the $(2n-6)$-dimensional components of $\mathcal{M}^*(Y, \operatorname{SL}_2(\C))$ (without the $(\Z_2)^{n-1}$ action). Moreover, the eigenvalues have multiplicities $\binom{n-1}{2}$. Together, the Milnor number is
\[
\mu(X) = \sum_{3 \leq n \leq k } e_n(a_1-1, \dots, a_k-1) \binom{n-1}{2},
\]
which gives the result by comparing with \eqref{eq:euler_char_moduli}.
\end{proof}

\begin{Example}
For $Y = \Sigma(2, 3, 5, 7)$, we have $e_3(1, 2, 4, 6)/4 = 23$ zero-dimensional components and $e_4(1, 2, 4, 6)/8 = 6$ two-dimensional components and
$
	\lambda^P_{\operatorname{SL}_2(\C)}(Y) = 23+6 \cdot 6 = 59$.
The corresponding Brieskorn--Hamm singularity $X$ has 92 eigenvalues of multiplicity 1 and 48 eigenvalues of multiplicity 3, so that $\mu(X) = 236 = 4 \cdot 59$.
\end{Example}

The above considerations suggest that this relation is not just a numerical coincidence and should be given an interpretation on the level of cohomology of the moduli space $\mathcal{M}^*(Y, \operatorname{SL}_2(\C))$ and of the Milnor fiber. Abouzaid and Manolescu defined the sheaf-theoretic $\operatorname{SL}_2(\C)$ Floer homology using Heegaard splitting, giving two Lagrangians in the character variety of the Heegaard surface. They considered certain perverse sheaf of vanishing cycles, which can be thought of as being associated with a function related to the Lagrangians.\footnote{In geometric and physics realizations of knot and 3-manifold invariants this function is known as the superpotential of a 3d theory \cite{ EGGKPSS22, EGGKPS22, FGSS12, GGP13a, GNSSS15,OV99, ERT22}; it defines a Landau--Ginzburg (LG) model that we discuss in more detail shortly.} The resulting cohomology is expected to be closely related to $\text{MTC} [Y, G_{\CC}]$ and to the Floer homology of Vafa--Witten theory on $Y$; see, e.g., \cite{GGP13b, GPV16, GSY22}. In particular, these connections indicate a direct relation to the sheaf counting (Vafa--Witten theory) on the complex surface $X$.

Further support for this relation comes from the analogy between the above proposition and the ``Casson invariant conjecture'' relating the $\lambda_{\operatorname{SU}(2)}$ with the one eight of the signature $\sigma$ of the Milnor fiber \cite{NeuWa90}. Since the latter has been verified for splice-quotients in \cite{NeOk08}, it would be interesting to explore a suitable analogue for $\operatorname{SL}_2(\C)$.

\subsection{Nonabelian noncommutative Hodge structure (NCNA)}

In this section, we give some questions and speculative interpretations of the relation between the moduli spaces of flat connections and the spectra of singularities. We will develop this further in separate paper.

Consider a germ of a singular affine hypersurface $X$ (more generally, a complete intersection) with an isolated singular point $0$, with the link $Y$ of $X$ being a homology sphere.\footnote{For rational homology spheres, one would have to pass to the universal abelian covering $X^{\rm ab}$, which makes the discussion more involved.} We have
$\pi_1(X {\setminus} 0) \cong \pi_1(Y)$,
so we can use the algebraic structure on $X {\setminus} 0$ to study moduli spaces of representations.

Following C.~Simpson \cite{Sim97} and T.~Mochizuki (see the survey of Sabbah \cite{Sab11}), we have a~real-analytic isomorphism of the Betti and Dolgachev moduli spaces
$H^1_B (Y) \cong M_{\rm Dol}(X {\setminus} 0)$.

\begin{Theorem}[Simpson, Mochizuki]
	$H^1_{\rm Dol}(X {\setminus} 0)$ carries a nonabelian mixed Hodge structure.
\end{Theorem}

As in the previous subsection, we consider the smoothing one-parameter family
$
	\{X_t \to X_0\}$,
with the monodromy operator (see Figure~\ref{fig:monodromy}). Associated with the family is the spectrum of singularity, denoted by $S_{\alpha}$, encoding the Hodge-theoretic information.

\begin{figure}[ht]
	\centering
	\begin{tikzpicture}
		\draw (-2, 0) arc (180:360:2 and 0.5);
		\draw (-2, 3) arc (180:-180:2 and 0.5);
		\draw (1, 1.5) node{$X_t$};
		\draw (-1, 1.5) node{$X_0$};
		\draw (-2, 0) --(-2, 3) (2, 0) --(2, 3);
		\draw[thick] (0, 0) to [out = 90, in = -20] (-0.5, 1.5) to [out = 20, in = -90] (0, 3);
		\draw[thick] (1, 0) to [out = 90, in = -90] (0.5, 1.5) to [out = 90, in = -90] (1, 3);
		\draw[-stealth] (1, 0) arc (-180:150:-1 and 0.25);
	\end{tikzpicture}
	\caption{The monodromy of an isolated surface singularity.}
	\label{fig:monodromy}
\end{figure}

Many of the classical aspects of Hodge theory have not been developed for mixed nonabelian Hodge structures. In particular, no Clemens--Schmid sequence is known for mixed nonabelian Hodge structures. We propose the nonabelian noncommutative Hodge structure (NCNA) as a possible approach to the Clemens--Schmid sequence for mixed nonabelian Hodge structures based on the theory of Landau--Ginzburg models. Our considerations suggest the following nonabelian noncommutative Hodge Structure, illustrated in Figure~\ref{fig:NCNA}.

\begin{figure}[ht]
	\centering
	\begin{tikzpicture}
		\draw (0, 0) node(b){$\CC$};
		\draw (0, 1.5) node(a){$\mch(X)$};
		\draw[-stealth] (a)-- (b) node[midway, left] {$\varphi$};
		\begin{scope}[shift = {(6, 0)}, scale = 1.2]
			\draw(-3, -0.3)--(1.8, -0.3)--(1.8, 2)--(-3, 2)--cycle;
			\draw (-3.5, 0.75) node{$\mch(X)$};
			\draw (-2.5, 0) node[scale = 0.8]{$\bullet$};
			\draw (-1.5, 0) node[scale = 0.8]{$\bullet$};
			\draw (-0.5, 0) node[scale = 0.8]{$\bullet$};
			\draw (-2.5, 1) node[scale = 0.8]{$\bullet$};
			\draw (-1.5, 1) node[scale = 0.8]{$\bullet$};
			\draw (-0.5, 1) node[scale = 0.8]{$\bullet$};

			\draw (0.5, 0) node[scale = 0.8]{$\bullet$};
			\draw (0.5, 1) ellipse (0.2 and 0.7) ;
			\draw (-1, 0) ellipse (1.9 and 0.3) ;

			\draw [decorate, decoration = {brace, amplitude = 8pt}, yshift = 15pt]
			(0.7, -0.9) -- (-2.7, -0.9)node [black, midway, yshift = -20pt] {$S_\alpha$};

			\draw (0.5, 1.4) -- (2, 1.4) node[right]{$H^1_{\rm Dol}(X {\setminus} 0)$};
			\draw (1.3, 0.5) node{$\mch_{\lambda_i, n}$};
		\end{scope}
	\end{tikzpicture}
	\caption{Nonabelian noncommutative structure combining the moduli spaces of representations with Landau--Ginzburg theory.}
	\label{fig:NCNA}
\end{figure}

The critical sets of the function $\varphi$ are the moduli spaces of connections, carrying the nonabelian Hodge structures. The proposals for constructing $\varphi$ using a reduction of the Chern--Simons functional were given, e.g., \cite{EGGKPSS22,EGGKPS22}. The noncommutative Hodge structure (NC) is encoded in the global structure of $\varphi$ and in the spectrum, which labels the components, and the multiplicity corresponds to the Euler characteristic of the components, as explained in the previous subsection. The noncommutative nonabelian Hodge structure $\mch$ can be thought of as a combination of these two phenomena. The properties of $\mch$ are expected to be combinations of A and B sides of properties \cite{KKP08}. We formulate these as questions.

\newcounter{questionsNCNA}

\begin{Question}\quad
	\begin{enumerate}\itemsep=0pt
		\item[(1)] Does $\mch$ behave as a B-side classical mixed Hodge structure?
		\item[(2)] Does $\mch$ behave as an A-side NC Hodge Structure?
		\item[(3)] Does $\mch$ have all properties of spectra, e.g., Thom--Sebastiani?
		\item[(4)] Are the stability conditions of the corresponding Fukaya--Seidel category determined by the bases of the Hitchin systems for the moduli spaces of representations?
		\item[(5)] Do the wall-crossing of phenomena for the above stability conditions determine the $\Zhat_{\sigma} (q)$ invariants associated with $Y$? This approach relates to the calculations by \cite{EGGKPSS22,EGGKPS22} and~\cite{KS1}.
		\setcounter{questionsNCNA}{\value{enumi}}
	\end{enumerate}
\end{Question}
Here are some additional questions about the abstract properties of proposed structures:
\begin{enumerate}
	\setcounter{enumi}{\value{questionsNCNA}}\itemsep=0pt
	\item[(1)] Consider an $L$-hyperplane section of $X$. Is it true that $\mch_{X \cap L} \subset \mch_X$?
	\item[(2)] Does $\mch$ have Thom--Sebastiani and semi-continuity properties?
	\item[(3)] Does the Hitchin base $B$ of the moduli space $M_{\rm Dol}(X {\setminus} 0)$ embed into the space of stability conditions for Fukaya--Seidel category associated to the LG model $\mch_X$?
	\[
		\operatorname{Stab}\mathop{FS}(\mch(X)) \supset B = \bigoplus\limits_{\lambda_i \in {\rm Spec}} H_{\lambda_i}\bigl({\rm Sym}^2(X {\setminus} 0)\bigr)?
	\]
\end{enumerate}

\begin{Example}
	$X\colon x^2 + y^3 + z^7 = 0 \supset$ Brieskorn homology sphere $\Sigma(2, 3, 7)$.
	\begin{center}
		\begin{tikzpicture}[scale = 1.2]
			\draw(-3, -0.5)--(1, -0.5)--(1, 2)--(-3, 2)--cycle;
			\draw (-3.5, 0.75) node{$\mch$};
			\draw (-2.5, 0) node[scale = 0.8]{$\bullet$};
			\draw (-1.5, 0) node[scale = 0.8]{$\bullet$};
			\draw (0, 0) node[scale = 0.8]{$\bullet$};
			\draw (0, 0) to[out = 80, in = 0, looseness = 1.2] (0, 1.5) to [out = 180, in = 100, looseness = 1.2] (0, 0) ;
			\draw (0, 1) -- (2, 1) node[right]{$\operatorname{SL}_2(\C)$};
		\end{tikzpicture}
	\end{center}
	There are three irreducible flat $\operatorname{SL}_2(\C)$-connections, corresponding to the three $(\Z_2)^2$ orbits of the spectrum of $X$ given by
$\frac{d_1}{2}+\frac{d_2}{3}+\frac{d_3}{7}
$ with $d_1 = 1$, $d_2 = 1, 2$ and $d_3 = 1, 2, 3, 4, 5, 6$.
	The Landau--Ginzburg model for the link of $X$ was proposed in \cite{EGGKPSS22}. As pointed out in loc.cit., the Landau--Ginzburg potential may have additional critical points that have interpretation (and play an important role) in curve counting and also in 3d-3d correspondence \cite{FGSS12, GGP13a}.
\end{Example}

Motivated by the previous subsection, we can see two approaches to building not only new NCNA Hodge structures but also two new ways to build invariants of 3-manifolds.

\emph{Approach $1$:} In \cite{AbMa20}, Abouzaid--Manolescu show that the moduli space of $\operatorname{SL}_2(\C)$ representations of $\pi_1(Y)$ (for any 3-manifold $Y$) is a derived critical locus, and they prove it by choosing a Heegaard splitting. The question is how the space and the function $f$ such that $\text{crit}(f) = X$ depend on the choice of Heegaard splitting.

In this derived situation, it is conceivable we can define a derived version of the classical spectrum.

\begin{Conjecture}
	The derived spectrum of $f$ is an invariant of $Y$, i.e., independent of Heegaard splitting of $Y$.
\end{Conjecture}

\emph{Approach $2$:} Use the approach of \cite{EGGKPSS22,EGGKPS22} to construct the potential. We must find a Lagrangian submanifold in a 3d Calabi--Yau manifold $T^* Y$ associated with $Y$; the disk counting then produces a potential $f$. Carrying this out might give a simpler proof of the above conjecture. Furthermore, the Fukaya--Seidel category for $f$ in this approach is expected to be related to the ``category of line operators'' $\text{MTC} [Y, G_{\CC}]$ in \cite{CCFGH18,GPV16}.

The spectra defined in the first and second approaches may be the same. It seems conceivable that the following conjecture holds.

\begin{Conjecture}
	The derived spectrum is the ``same'' for all potentials $f$ with the same derived critical locus $X$.
\end{Conjecture}

The above considerations suggest three invariants of derived singularity theory with increasing order of complexity.
\begin{enumerate}\itemsep=0pt
	\item[(1)] Hypercohomologies of the perverse sheaf of vanishing cycles $F$. When $F$ comes from the Heegaard splitting of a 3-manifold $Y$, the above hypercohomologies form an invariant of~$Y$~\cite{AbMa20}.

	\item[(2)] We can enhance the above hypercohomologies with a mixed Hodge structure. Combined with the monodromy, it leads to the \emph{derived spectrum}. In the case of the derived spectrum connected with a moduli space of $\operatorname{SL}_2(\C)$ representations, it defines invariants of this moduli space and an invariant of the 3-manifold $Y$.

	\item[(3)] We can bring additional data to the above potential $f$ --- a divisor $D$ along which $f$ has a~log behavior. This allows an additional spectral grading from the number of times a~path coming from one component of the critical set goes around $D$ and ending on another component of the critical set.
\end{enumerate}

Changing the divisor $D$ leads to different filtrations and functors among different Fukaya--Seidel categories.
In such a way, we get an enhancement of classical singularity theory. In the case of a 3-manifold, all of the above data of singularity theory (commutative spectrum, derived spectrum, and spectral grading) are recorded by the $q$-series $\Zhat_{\sigma} (Y;q)$. These make $\Zhat_{\sigma} (Y;q)$ an interesting starting point for producing invariants of 4-manifolds via categorification.

Some of these ideas have appeared before in the singularity theory and category theory --- see, e.g., \cite{FichYin18,Kuz09}. This will require a new theory of spectra in the case of shifted symplectic structures.

\subsection{Relation to vertex algebras and invariants of complex surfaces}\label{sec:ADEVOA}
The nonconventional modular properties of $\Zhat_{\sigma} (q)$-invariants perfectly fit those exhibited by characters of a logarithmic vertex algebra. In fact, the relation of the form
$
	\Zhat_{\sigma} (q) = \chi_{\sigma} (q)
$
is one of the predictions of the so-called ``3d-3d correspondence'' in physics, and can be viewed as its mathematical incarnation. Here, $\chi_{\sigma} (q)$ is a character of a log-VOA labeled by 3-manifold $Y$ (and a choice of the root system, which is implicit throughout the paper).

For a plumbed 3-manifold $Y$, a specific choice of a plumbing graph corresponds to a choice of a 4-manifold $X$, with $Y = \partial X$ as its boundary ($X$ is given by 4d Dehn surgery). If the plumbing is negative definite, this 4-manifold can be equipped with a complex structure (as a resolution of a normal surface singularity) and we can study analytical invariants of the complex surfaces and compare them to topological invariants of the underlying 4-manifold and to invariants of~$Y$. Since various $q$-series invariants of $X$ were available for quite some time, it is natural to ask how they compare with $\Zhat_{\sigma} (q)$ and $Z_0 (q)$. The development of $\Zhat_{\sigma} (Y;q)$ was, in fact, largely motivated \cite{GPV16} by the connection to the Vafa--Witten and Donaldson--Thomas invariants of $X$ which share many similarities with $\Zhat_{\sigma} (Y;q)$. For example, both require a choice of an additional structure whose cutting-and-gluing is described by a `decorated' version of TQFT, and both are related to characters of VOAs labeled by the corresponding 3-manifolds or 4-manifolds.

To keep this discussion more concrete and to explore these relations more deeply, let us consider a rather special class of spherical 3-manifolds of $ADE$ type, i.e., $Y = S^3 / \Gamma$, where $\Gamma$ is a discrete subgroup of $\operatorname{SU}(2)$. In type $A$, $\Gamma = \ZZ_p$ is simply a cyclic group and $Y = L(p, 1)$ is a Lens space. Similarly, for $E_8$ we get a Poincar\'e sphere, and we refer to other spherical 3-manifolds in this family by their type, e.g., as ``$E_6$ manifold'' and so on.

\begin{Theorem}[following \cite{Nak94}]\label{theorem:ALE}
	For a discrete subgroup $\Gamma \subset \operatorname{SU}(2)$ of ADE type, let $X_{\Gamma}$ be the corresponding ALE space given by a resolution of the $\CC^2 / \Gamma$ singularity. Then, for any
\[\rho \in \operatorname{Hom} ( \Gamma , {\rm GL}(N) )
\]
	the generating series of Vafa--Witten $($and Donaldson--Thomas$)$ invariants counting rank-$N$ sheaves on $X_{\Gamma}$ is a character of the affine Lie algebra $\widehat{\mathfrak{g}}$ at level $N$
	\[
		Z_{VW} (X_{\Gamma}; q, \rho) =
		Z_{DT} (X_{\Gamma}; q, \rho) =
		\chi_{\rho}^{\widehat{\mathfrak{g}}_N} (q),
	\]
	where $\mathfrak{g}$ is related to $\Gamma$ by McKay correspondence.
\end{Theorem}

Note that all representations $\rho\colon \Gamma \to {\rm GL}(N)$ can be conjugated to ${\rm U}(N)$ representations. In particular, this means that for spherical manifolds $S^3 / \Gamma$ the problem of enumerating ${\rm U}(N)$ flat connections is equivalent to the problem of classifying ${\rm GL}(N)$ flat connections.\footnote{This simple fact plays an important role in relating infinite-dimensional algebras and modular data associated to 4-manifolds and to 3-manifolds \cite{FG18}; in particular, for spherical 3-manifolds the set $\{ \rho \}$ labels simple objects in $\text{MTC} [Y]$ \cite{CCFGH18}.} For example, when $\Gamma$ is of type $A$, i.e., a cyclic group $\ZZ_p$, we can think of $\rho$ as a Young tableau that can fit in a rectangle (of the size determined by $N$ and the order of $\Gamma$). The number of such $\rho$ is equal to \smash{$\frac{(N+p-1)!}{(p-1)! N!}$}. We can also use the relation \smash{${\rm U}(N) = \frac{{\rm U}(1) \times \operatorname{SU}(N)}{\ZZ_N}$} to enumerate $\operatorname{SU}(N)$ or $\operatorname{SL}(N)$ representations
\smash{$\rho\colon \Gamma \to \operatorname{SL}(N)$}.
Their number, equal to \smash{$\frac{(N+p-1)!}{(N-1)! p!}$}, is the number of representations of \smash{$\widehat{\mathfrak{su} (N)}_p$}.

The analogous statements for 3-manifold $q$-series invariants can be stated in terms of \smash{$\Zhat_{\sigma} (Y; q)$}. We will illustrate it on the $E_6$ manifold, continuing Example~\ref{ex:E6}. Other spherical manifolds can be treated similarly.
\begin{Example}
	From the expressions in Example~\ref{ex:E6}, we can write $\Zhat_0$ and $\Zhat_1$ in terms of false theta-functions
	\begin{align}
		\Zhat_0 (E_6) & = \frac{2 q^{-1}}{(q)_{\infty}} - \frac{q^{- \frac{25}{24}}}{(q)_{\infty}} \bigl( \tilde \Psi^{(1)}_6 + \tilde \Psi^{(5)}_6 \bigr) = \frac{q^{-1}}{(q)_{\infty}} \bigl(
		1 - q + q^2 + q^5 - q^7
		+ \cdots \bigr),
		\nonumber\\
		\Zhat_1 (E_6) & = \frac{q^{ - \frac{2}{3} - \frac{3}{8}}}{(q)_{\infty}} \tilde \Psi^{(3)}_6 = \frac{q^{-2/3}}{(q)_{\infty}} \bigl(
		-1 + q^3 - q^9 + q^{18} - q^{30} + q^{45}
		+ \cdots \bigr) ,\label{ZhatE6}
	\end{align}
	where we restored the universal factor $\frac{1}{(q)_{\infty}}$ present for all 3-manifolds\footnote{Due to its universal nature, this factor is often omitted for simplicity, but it is important in matching partition functions of 3d physical theory \cite{GPPV20} and in comparing with characters of logarithmic vertex algebras~\cite{CCFGH18, CCFFGHP22, Sug22}. If we were to study higher rank invariants for more general root systems, there would be several such factors \cite{Pa20a}.} and used the standard notation for the false theta-functions
	\begin{gather}
		\tilde \Psi^{(a)}_p (q) : = \sum_{n = 0}^\infty \psi^{(a)}_{2p}(n) q^{\frac{n^2}{4p}} \in q^\frac{a^2}{4p} \Z[[q]],
		\qquad
		\psi^{(a)}_{2p}(n) =
		\begin{cases}
			\pm 1, & n\equiv \pm a~\mod~ 2p , \\
			0, & \text{otherwise}.
		\end{cases} \label{falsetheta}
	\end{gather}
	In particular, we have
$\Delta_\sigma = h_\sigma - \frac{c}{24} \in \bigl\{ -1 , - \tfrac{2}{3} \bigr\}$,
	which already gives some information about the central charge $c$ of the corresponding log-VOA and conformal weights $h_\sigma$ of its modules. Let~${a, b \in \mathbb{Z}_{\geq 1}}$, $\gcd(a, b) = 1$. For $r', s' \in \mathbb{Z}$, set
\[
\beta_{r',s'} := \left( -\frac{r'-1}{a} + \frac{s'-1}{b} \right) \sqrt{ab}\varpi,
\]
where $\varpi = \frac{1}{2}\alpha$ is the fundamental weight of $\mathfrak{g} = \mathfrak{sl}_2$. In particular, for $n, n_1, n_2 \in \mathbb{Z}$, $n_1 + n_2 = n$ and $1 \leq s \leq a$, $1 \leq r \leq b$, we have
\[
\beta(-n, s, r) := \left( -n - \frac{s-1}{a} + \frac{r-1}{b} \right) \sqrt{ab}\varpi = \beta_{2-r+bn_1, 2-s-an_2}
\]
and let us denote by $\pi_\beta$ the Fock space with weight $\beta$. For the conformal vector
\[
\omega := \frac{1}{2}\alpha_{(-1)} \varpi + Q_0 \varpi_{(-2)} |0\rangle \in \pi_0, \qquad Q_0 := \sqrt{ab}\left(\frac{1}{a} - \frac{1}{b}\right) = \sqrt{\frac{b}{a}} - \sqrt{\frac{a}{b}},
\]
the central charge $c$ and the conformal weight $h_{r',s'}$ of $|\beta_{r',s'}\rangle \in \pi_{\beta_{r',s'}}$ is given by
\begin{gather*}
c = 1 - 12|Q_0 \varpi|^2 = 1 - \frac{6(a-b)^2}{ab} = 1 - 6\left( \frac{a}{b} - 2 + \frac{b}{a} \right),
\\
h_{r',s'} = \frac{1}{2} |\beta_{r',s'} - Q_0 \varpi|^2 + \frac{c-1}{24} = \frac{(ar'-bs')^2-(a-b)^2}{4ab}
\\
\phantom{h_{r',s'} }{}
= - \frac{r'^2-1}{4}\frac{b}{a} - \frac{r's'-1}{2}+ \frac{s'^2-1}{4}\frac{a}{b},
\end{gather*}
Therefore, the character $\chi_M := \text{tr}_M q^{L_0 - \frac{c}{24}}$ of $\pi_{\beta_{r',s'}}$ is
\[
\chi_{\pi_{\beta_{r',s'}}} = \frac{1}{(q)_\infty} q^{h_{r',s'} - \frac{c}{24}} = \frac{1}{\eta(q)} q^{\frac{1}{2}|\beta_{r',s'} - Q_0 \varpi|^2}.
\]
In particular,
\[
\chi_{\pi_{\beta(-n;s,r)}} = \frac{1}{\eta(q)} q^{\frac{1}{2}|\sqrt{ab}(-n-\frac{s}{a}+\frac{r}{b})\varpi|^2} = \frac{1}{\eta(q)} q^{\frac{1}{4 a b}(-nab-bs+ar)^2}.
\]
For the case $(a,b) = (1,p)$ with $p \geq 2$ (so that $s = 1$ and $1 \leq r \leq p$), we have
\[
\chi_{\pi_{\beta(-n;1,r)}} = \frac{1}{\eta(q)} q^{\frac{1}{4p}(-(n+1)p+r)^2}.
\]
By the Felder complex, the character of the irreducible module $M(0; 1, r) \subset \pi_{\beta(0;1,r)}$ of the $(1, p)$ singlet log-VOA $M(0; 1, 1) \subset \pi_{\beta(0;1,1)}$ is given by
\[
\chi_{M(0;1,r)} = \sum_{m \geq 0} \chi_{\pi_{\beta(-2m;1,r)}} - \chi_{\pi_{\beta(-(2m+1);1,p-r)}} = \frac{1}{\eta(q)} \sum_{m \geq 0} q^{\frac{1}{4p}(-(2m+1)p \pm r)^2} = \tilde{\Psi}_p^{(p-r)}.
\]
For example, from \eqref{ZhatE6} we read off the value of $p = 6$, which gives $c = 13 - 6p - 6p^{-1} = -24$. This, in turn, implies that \smash{$h_\sigma = \Delta_\sigma + \frac{c}{24} = \frac{(p-s)^2}{4p} + \frac{c-1}{24} = \bigl\{ - 1, - \frac{2}{3} \bigr\}$} for the $E_6$ manifold.
	Therefore, we conclude that for the $E_6$ manifold, $\Zhat_{\sigma} (Y; q)$ are characters of a logarithmic $(1, 6)$ singlet model.
\end{Example}

Similarly, one can prove the following analogue of Theorem~\ref{theorem:ALE}.
\begin{Theorem}[following \cite{CCFFGHP22}]
	For the family of spherical $3$-manifolds, $Y = S^3 / \Gamma$, the non-perturbative $\operatorname{SL}_2(\C)$ invariants $\Zhat_{\sigma} (Y; q)$ are equal to characters of the following $($logarithmic$)$ vertex algebras:
	\begin{itemize}\itemsep=0pt
\item In type A, $\Zhat_{\sigma} (Y; q)$ is a character of an ordinary $($non-logarithmic$)$ vertex algebra, namely the character of a Feigin--Fuchs module of the Virasoro algebra $($or, equivalently, that of a~${c = 1}$ ``free boson'' VOA$)$.

		\item In type D$_n$, $\Zhat_{\sigma} (Y; q)$ is a character of a logarithmic $(1, p)$ singlet model, with $p = n-3$.

		\item In type E, $\Zhat_{\sigma} (Y; q)$ is a character of a logarithmic $(1, p)$ singlet model, with $p = 6$, $12$, and~$30$ in the case of $E_6$, $E_7$, and $E_8$, respectively.

	\end{itemize}
\end{Theorem}

\begin{Remark}
	In type $A$, the assertion follows from the explicit expression (see, e.g., \cite{GPV16})
	\[
		\Zhat_{\sigma} (Y; q) = \frac{1}{(q)_{\infty}} q^{\Delta_\sigma}.
	\]
	The only interesting aspect of this case are the values of $\Delta_\sigma = h_\sigma - \frac{c}{24}$, which are analogs of ``correction terms'' in the Heegaard Floer theory and whose computation requires some care~\mbox{\cite{AJK21, GPP20}}.
\end{Remark}
We also note that Theorem~\ref{thm:seifert} gives an explicit way of how to write $\Zhat_{\sigma} (q)$ invariants in terms of linear combinations of derivatives of false theta functions. We hope these expressions can help explore further relations between $\Zhat_{\sigma} (q)$ invariants and VOA characters and construct new logarithmic vertex algebras.

\subsection*{Acknowledgements}

We would like to thank Denis Auroux, Martin \v{C}ech, Shimal Harichurn, Mrunmay Jagadale, Maxim Kontsevich, Slava Krushkal, Andr\'as N\'emethi, Sunghyuk Park, Pavel Putrov and Nikolai Saveliev for helpful conversations and Zuzana Urbanov\'a for the help with computer experiments. We thank the anonymous referees for valuable suggestions and comments.
S.~Gukov was partially
supported by NSF grant DMS-1664227, DOE grant DE-SC0011632. S.~Gukov and J.~Svoboda were partially suported by the Simons Grant {\it New structures in low-dimensional topology}.
L.~Katzarkov was supported by NSF Grant theory of Atoms, Simons investigators grant
HMF Simons Foundation, grant SFI-MPS-T-Institutes-00007697, and the Ministry of Education and Science of the Republic of Bulgaria, grant DO1-239/10.12.2024.

\pdfbookmark[1]{References}{ref}
\LastPageEnding

\end{document}